\documentclass[11pt]{amsart}
\usepackage{amsmath,amsthm,amssymb,amsfonts,amscd,fullpage}

\newcommand \Tor {\ensuremath{\mathrm{Tor}}}
\newcommand \lex {\ensuremath{\mathrm{lex}}}
\newcommand \soc {\ensuremath{\mathrm{soc}}}
\newcommand \ini {\ensuremath{\mathrm{in}}}
\newcommand \Gin {\ensuremath{\mathrm{Gin}}}

\newcommand \bh {\mathbf H}

\newcommand \kxn[1] {k[x_1,\ldots,x_{#1}]}

\theoremstyle{plain} 
\newtheorem{thm}{Theorem}[section]
\newtheorem{prop}[thm]{Proposition}

\newtheorem{lem}[thm]{Lemma}
\newtheorem{cor}[thm]{Corollary}
\newtheorem{pro}[thm]{Proposition}

\theoremstyle{definition}

\newtheorem{remk}[thm]{Remark}

\newtheorem{ques}[thm]{Question}
\newtheorem{exmp}[thm]{Example}

\newtheorem{rem}[thm]{Remark}
\newtheorem{que}[thm]{Question}

\def\H{\text{\bf H}}
\def\P{{\mathbb P}}

\def\H{{\bf H}}

\def\ni{\noindent}

\def\ra{\rightarrow}

\def\ds{\displaystyle}

\def\FF{{\mathcal F}}

\def\GG{{\mathcal G}}

\numberwithin{equation}{section}

\begin{document}

\title[Generic Initial Ideals and Graded Artinian Level Algebras]{Generic Initial Ideals And Graded Artinian Level Algebras Not Having The Weak-Lefschetz Property}

\author[J.M. Ahn]{Jea-Man Ahn}
\address{Korea Institute for Advanced Study,  Seoul, Korea, 130-722}
\email{ajman@kias.re.kr}
\author[Y.S. Shin]{Yong Su Shin}
\address{Department of Mathematics,  Sungshin Women's University, Seoul,   Korea, 136-742}
\email{ysshin@sungshin.ac.kr}
\begin{abstract} We find a sufficient condition that $\H$ is not level based on a reduction number. In particular, we prove that a  graded Artinian algebra of codimension $3$ with Hilbert function $\H=(h_0,h_1,\dots, h_{d-1}>h_d=h_{d+1})$ cannot be level if $h_d\le 2d+3$, and that there exists a level O-sequence of codimension $3$ of type $\H$ for $h_d \ge 2d+k$ for $k\ge 4$. Furthermore, we show that $\H$ is not level if $\beta_{1,d+2}(I^{\rm lex})=\beta_{2,d+2}(I^{\rm lex})$, and also prove that any codimension $3$ Artinian graded algebra $A=R/I$ cannot be level if $\beta_{1,d+2}(\Gin(I))=\beta_{2,d+2}(\Gin(I))$. In this case, the Hilbert function of $A$ does not have to satisfy the condition $h_{d-1}>h_d=h_{d+1}$.

Moreover, we show that every codimension $n$ graded Artinian level algebra having the Weak-Lefschetz Property has the strictly unimodal Hilbert function having a growth condition on $(h_{d-1}-h_{d}) \le  (n-1)(h_d-h_{d+1})$ for every $d > \theta$ where
$$
h_0<h_1<\cdots<h_\alpha=\cdots=h_{\theta}>\cdots>h_{s-1}>h_s.
$$
In particular, we find that if $A$ is of codimension $3$, then $(h_{d-1}-h_{d}) <  2(h_d-h_{d+1})$ for every $\theta< d <s$ and $h_{s-1}\le 3 h_s$, and prove that if $A$ is a codimension $3$ Artinian algebra with an $h$-vector $(1,3,h_2,\dots,h_s)$ such that
$$
h_{d-1}-h_d=2(h_d-h_{d+1})>0 \quad \text{ and } \quad
\soc(A)_{d-1}=0
$$
for some $r_1(A)<d<s$, then $(I_{\le d+1})$ is $(d+1)$-regular and $\dim_k\soc(A)_d=h_d-h_{d+1}$.
\end{abstract}

\subjclass[2000]{Primary:13D40; Secondary:14M10}
\keywords{Level algebra, Gorenstein algebra, Betti number, Hilbert function, Weak Lefschetz Property, Generic Initial Ideal, Strictly Unimodal O-sequence.}

\maketitle

\section{Introduction}

Let $R=k[x_1,\dots,x_n]$ be an $n$-variable polynomial ring over an infinite field with characteristic $0$. In this article, we shall study Artinian quotients $A=R/I$ of $R$ where $I$ is a homogeneous ideal of $R$. These rings are often referred to as standard graded algebras. Since $R = \oplus_{i=0}^\infty R_i$ ($R_i$: the vector space of dimension $\binom{i+(n-1)}{n-1}$ generated by all the monomials in $R$ having degree $i$) and $I = \oplus_{i=0}^\infty I_i$, we get that
$$
A = R/I = \oplus_{i=0}^\infty (R_i/I_i) = \oplus _{i=0}^\infty A_i
$$
is a graded ring.  The numerical function
$$
\H_A(t):= \dim _k A_t = \dim _kR_t - \dim _k I_t
$$
is called the {\em Hilbert function} of  the ring $A$.

Given an O-sequence $\H=(h_0,h_1,\dots)$, we define the {\em first difference} of $\H$ as
$$
\Delta \H=(h_0,h_1-h_0,h_2-h_1,h_3-h_2,\dots).
$$

If $I$ is a homogeneous ideal of $R$ of height $n$, then $A=R/I$ is an {\em Artinian $k$-algebra}, and hence $\dim_k A<\infty$. We associate to the graded algebra $A$ a vector of nonnegative integers which is an $(s+1)$-tuple, called the {\em$h$-vector} of $A$ and denoted by
$$
h(A)=(h_0,h_1,\dots,h_s),
$$
where $h_i=\dim_k A_i$. Thus we can write $A=k\oplus A_1\oplus \cdots \oplus A_s$ where $A_s\ne 0$. We call $s$ the {\em socle degree} of $A$.  The {\em socle} of $A$ is defined by the annihilator of the maximal homogeneous ideal, namely
$$
{\rm ann}_A(m):=\{a\in A \mid am=0\} \qquad \text{where} \qquad m=\sum_{i=1}^s A_i.
$$
Moreover, an $h$-vector $(h_0,h_1,\dots,h_s)$ is called
$$
\begin{array}{llllllllllllllll}
\text{{\em unimodal} if } & h_0\le \cdots \le h_t=\cdots=h_\ell \ge \cdots \ge h_s,\\
\text{{\em strictly unimodal} if } & h_0< \cdots < h_t=\cdots=h_\ell > \cdots > h_s.
\end{array}
$$

A graded Artinian $k$-algebra $A=\bigoplus_{i=0}^s A_i$ ($A_s\ne 0$) is said to have  the {\em  Weak Lefschetz Property} (WLP for short) if there is an element $L\in A_1$ such that the linear transformations
$$
A_i \overset{\times L}{\rightarrow} A_{i+1}, \quad 1\le i \le s-1,
$$
which is defined by multiplication by $L$, are either injective or surjective. This implies that the linear transformations have maximal ranks for every $i$. In this case, we call $L$ a {\em Lefschetz element}.

A monomial ideal $I$ in $R$ is {\em stable} if the monomial
$$
\frac{x_jw}{x_{m(w)}}
$$ belongs to $I$ for every monomial $w\in I$ and $j<m(w)$ where
$$
m(u):=\max\{ j \mid a_j>0\}
$$
for $u=x_1^{a_1}\cdots x_n^{a_n}$. Let $S$ be a subset of all monomials in $R=\bigoplus_{i\ge 0} R_i$ of degree $i$. We call $S$ a {\em Boreal fixed set} if
$$
u=x_1^{a_1}\cdots x_n^{a_n}\in S,\ a_i>0,\ \quad \text{implies} \quad
\frac{x_i u}{x_j} \in S
$$
for every $1\le i\le j\le n$.

A monomial ideal $I$ of $R$ is called a {\em Borel fixed ideal or strongly stable ideal} if the set of all monomials in $I_i$ is a Borel set for every $i$. There are two Borel fixed monomial ideals canonically attached to a homogeneous ideal $I$ of $R$: the generic initial ideal $\Gin(I)$ with respect to the reverse lexicographic order and the
lex-segment ideal $I^{\rm lex}$.  The ideal $I^{\rm lex}$ is defined as follows. For the vector space $I_d$ of forms of degree $d$ in $I$, one defines $(I^{\rm lex})_d$ to be the vector space generated by largest, in the lexicographical order, $\dim_k(I_d)$ monomials of degree $d$. By construction, $I^{\rm lex}$ is a strongly stable ideal and it only depends on the Hilbert function of $I$.

In case of the generic initial ideal, it has been proved that generic initial ideals are Borel fixed in characteristic zero by Galligo~\cite{Ga}, and then generalized by Bayer and Stillman to every characteristic \cite{BS}.

In \cite{AM}, they gave some geometric results using generic initial ideals for the degree reverse lexicographic order, which improved a well known result of Bigatti, Geramita, and Migliore concerning geometric consequences of maximal growth of the Hilbert function of the Artinian reduction of a set of points in \cite{BGM}.
In~\cite{GHMS}, they gave a homological reinterpretation of   a level Artinian algebra and explained the combinatorial notion of Cancellation of Betti numbers of the minimal free resolution of the lex-segment ideal associated to a given homogeneous ideal. We shall explain the new result when we carry out the analogous result using the generic initial ideal instead of the lex-segment ideal. We find some new results on the maximal growth of the difference of Hilbert function in degree $d$ larger than the reduction number $r_1(A)$ if there is no socle element in degree $d-1$ using some recent result in \cite{AM}. As an application, we  give the condition if some O-sequence is  ``either level or non-level sequences¡± of Artinian graded algebras with the WLP.

\medskip

Let $\FF$ be the graded minimal resolution of $R/I$, i.e.,
$$
\begin{array}{ccccccccccccccccccccccccccccccccccccccccccc}
\FF: & 0 & \ra & \FF_n & \ra & \FF_{n-1} & \ra & \cdots
     & \ra & \FF_1 & \ra & R & \ra & R/I & \ra & 0.
\end{array}
$$
We can write
$$
\FF_i
 =  \bigoplus^{\gamma_{i}}_{j=1} R^{\beta_{ij}}(-\alpha_{ij})
$$
where $\alpha_{i1}< \alpha_{i2}< \cdots < \alpha_{i\gamma_{i}}$. The numbers $\alpha_{ij}$ are called the {\em shifts} associated to $R/I$, and the numbers $\beta_{ij}$ are called the {\em graded Betti numbers} of $R/I$. For $I$ as above, the {\em Betti diagram} of $R/I$ is a useful
device to encode the graded Betti numbers of $R/I$ (and hence of $I$).
It is constructed as follows:
$$
\bordermatrix{&  & 0 & 1 & \cdots & n-1 \cr 0 & 1 & 0 & 0 & \cdots & 0 \cr
1 & 0 & * & * & \cdots & * \cr
\vdots& \vdots & \vdots & \vdots & \vdots & \vdots \cr t & 0 &
\beta_{0,t+1} & \beta_{1,t+2} & * & \beta_{n-1, t+n} \cr
\vdots & \vdots & \vdots & \vdots & \vdots & \vdots  \cr d-2 & 0
&\beta_{0,d-1} & \beta_{1,d} & * & \beta_{n-1,d-2 + n} \cr d-1 & 0 &
\beta_{0,d} & \beta_{1,d+1} & * & \beta_{n-1, d-1+n} \cr d & 0 &
\beta_{0,d+1} & \beta_{1,d+2} & * &\beta_{n-1,d+n} \cr
\vdots & \vdots &\vdots & \vdots & \vdots &  }
$$
When we need to emphasize the ideal $I$, we shall use $\beta_{i,j}(I)$ for $\beta_{i,j}$.

Now, we recall that if the last free module of the minimal free resolution of a graded ring $A$ with Hilbert function $\H$ is of the form $\FF_n=R^\beta(-s)$ for some $s>0$, then   Hilbert function $\H$ and a graded ring $A$ are called {\em level}. For a special case, if $\beta=1$, then we call a graded Artinian algebra $A$ {\em Gorenstein}. In~\cite{St}, Stanley proved that any graded Artinian Gorenstein algebra of codimension $3$ is unimodal. In fact, he proved a stronger result than unimodality using the structure theorem of Buchsbaum and Eisenbud for the Gorenstein algebra of codimension $3$ in \cite{BE}. Since then, the graded Artinian Gorenstein algebras of codimension $3$ have been much studied (see \cite{D}, \cite{GHMS}, \cite{GHS:2},  \cite{Ha:1}, \cite{Ha:2}, \cite{M}, \cite{MN-1}, \cite{Sh:1}, \cite{Za1}). In~\cite{BI}, Bernstein and Iarrobino showed how to construct non-unimodal graded Artinian Gorenstein  algebras of codimension higher than or equal to $5$. Moreover, in \cite{BL}, Boij and Laksov showed another method on how to construct the same graded Artinian Gorenstein  algebras. Unfortunately, it has been unknown if there exists a graded non-unimodal Gorenstein algebra of codimension $4$. For unimodal Artinian Gorenstein algebras of codimension $4$, how to construct  some of them using the link-sum method has been shown in \cite{Sh:1}. It has been also shown in \cite{GHS:2} and \cite{Ha:1} how to obtain some of unimodal Artinian Gorenstein algebras of any codimension $n\ (\ge 3)$. An SI-sequence is a finite sequence of positive integers which is symmetric, unimodal and satisfies a certain growth condition.
In \cite{MN-1}, Migliore and Nagel showed how to construct a reduced, arithmetically Gorenstein configuration $G$ of linear varieties of arbitrary dimension whose Artinian reduction has the given SI-sequence as Hilbert function and has the Weak Lefschetz Property.
For graded Artinian level algebras, it has  been recently studied (see \cite{BI},   \cite{BG}, \cite{BL}, \cite{Cho-Iarrobino}, \cite{GHMS}, \cite{GHS:4}, \cite{M}, \cite{Za1}, \cite{Za2}).
In~\cite{GHMS}, they proved the following result.  Let
\begin{equation} \label{EQ:011}
\begin{array}{llllllllllllllllllllll}
\H & : & h_0 & h_1 & \cdots & h_{d-1} & h_d & h_d & \cdots
\end{array}
\end{equation}
with $h_{d-1}>h_d$. If  $h_d\le d+1$ with any codimension $h_1$, then $\H$ is {\em not} {level}.

In \cite{Za1}, F. Zanello constructed a non-unimodal level O-sequence of codimension $3$ as follows:
$$
\H=(h_0,h_1,\dots,h_d,t,t,t+1,t,t,\dots,t+1,t,t)
$$
where the sequence $t,t,t+1$ can be repeated as many times as we want. Thus there exists a graded Artinian level algebra of codimension $3$ of type in equation~\eqref{EQ:011} which does not have  the WLP.

\medskip

In Section 2, preliminary results and notations on lex-segment ideals and generic initial ideals are introduced. In Section 3, we show that any codimension $n$ graded Artinian level algebra $A$ having the  WLP has the Hilbert function which is strictly unimodal (see Theorem~\ref{T:307}). In particular, we prove that if $A$ has the Hilbert function such that
$$
h_0<h_1<\cdots<h_{r_1(A)}=\cdots=h_\theta>\cdots>h_{s-1}>h_s,
$$
then $h_{d-1}-h_d\le (n-1)(h_d-h_{d+1})$ for every $\theta<d \le s$ (see Theorem~\ref{T:307}). Furthermore, we show that if $A$ is of codimension $3$, then $h_{d-1}-h_d < 2(h_d-h_{d+1})$ for every $\theta<d <s$ and $h_{s-1}\le 3h_s$ (see Theorem~\ref{T:303}). We also prove that if $A$ is a codimension $3$ Artinian graded algebra with socle degree $s$ and
$$
\beta_{1,d+2}(\Gin(I))=\beta_{2,d+2}(\Gin(I))>0
$$
for some $d<s$, then $A$ cannot be level (see Theorem~\ref{T:301-1}). Moreover, if $A=R/I$ is a codimension $3$ Artinian graded algebra with an $h$-vector $(1,3,h_2,\dots,h_s)$ such that $h_{d-1}-h_d=2(h_d-h_{d+1})>0$ for some $r_1(A)<d<s$ and $\soc(A)_{d-1}=0$, then $(I_{\le d+1})$ is $(d+1)$-regular and $\dim_k \soc(A)_d=h_d-h_{d+1}$ (see Theorem~\ref{T:302}).

\medskip

One of the main topics of this paper is to study O-sequences of type in equation~\eqref{EQ:011} and find an answer to the following question.

\begin{que} \label{Q:111}
Let $\H$ be as in equation~\eqref{EQ:011} with $h_1=3$. What is the
minimum value for $h_d$ when $\H$ is level?
\end{que}

Finally in Section~\ref{Sec:004}, we show that if $R/I$ is a graded Artinian algebra of codimension $3$ having  Hilbert function $\H$ in equation~\eqref{EQ:011} and $\beta_{1,d+2}(I^{\rm lex})=\beta_{2,d+2}(I^{\rm lex})$, then $R/I$ is {\em\bfseries not} level, i.e., $\H$ cannot be level (see Theorem~\ref{T:407}).  Furthermore, we prove that any O-sequence $\H$ of codimension $3$ in equation~\eqref{EQ:011} cannot be level when $h_d\le 2d+3$ and there exists  a level O-sequence of codimension $3$ of type in equation~\eqref{EQ:011} having $h_d\ge 2d+k$ for every $k\ge 4$ (see Theorem~\ref{T:403}, Proposition~\ref{P:411} and Remark~\ref{R:412}), which is a complete answer to Question~\ref{Q:111}.

\medskip

A computer program CoCoA was used for all examples in this article.

\section{Some Preliminary Results}


In this section, we introduce some preliminary results and notations on lex-segment ideals and generic initial ideals.

\begin{thm}[\cite{AM}, \cite{BS}, \cite{G}]\label{Generic initial ideal of hyperplane section}
Let $L$ be a general linear form and let
$J=(I+(L))/(L)$ be considered as a homogeneous ideal of
$S=k[x_1,\ldots,x_{n-1}]$. Then
$$
\Gin(J)=\big(\Gin(I)+(x_n)\big)/(x_n).
$$
\end{thm}

For a homogeneous ideal $I\subset R$ there exists a flat family
of ideals $I_t$ with $I_0=\ini(I)$ (the initial ideal of $I$) and $I_t$ canonically
isomorphic to $I$ for all $t\neq 0$ (this implies that $\ini(I)$
has the same Hilbert function as the one of $I$). Using this
result, we get the following Theorem:

\begin{thm}[The Cancelation Principle, \cite{AM}, \cite{G}]\label{Cancellation Principle}
For any homogeneous ideal $I$ and any
$i$ and $d$, there is a complex of $k\cong R/m-$modules
$V_{\bullet}^d$ such that
$$
\begin{array}{rclllllllll}
V_{i}^d
& \cong & \Tor_{i}^R(\ini(I),k)_d \\[1ex]
H_i(V_{\bullet}^d)
& \cong & \Tor_{i}^R(I,k)_d.
\end{array}
$$
\end{thm}

\begin{remk}\label{remark of CP}
One way to paraphrase this Theorem is to say that the minimal
free resolution of $I$ is obtained from that of $\ini(I)$, the {\em initial ideal} of $I$, by
canceling some adjacent terms of the same degree.
\end{remk}

\begin{thm}[Eliahou-Kervaire, \cite{EK}]\label{EK}
Let $I$ be a stable monomial ideal of $R$. Denote by $\mathcal
G(I)$ the set of minimal (monomial) generators of $I$ and
$\mathcal G(I)_d$ the elements of $\mathcal G(I)$ having degree
$d$. Then
$$
\beta_{q,i}(I)=\sum_{T\in \mathcal
G(I)_{i-q}}\binom{m(T)-1}{q}.
$$
\end{thm}

This theorem gives all the graded Betti numbers of the lex-segment ideal and the generic initial ideal just from an intimate knowledge of the generators of that ideal. Since the minimal free resolution of the ideal of a $k$-configuration in $\P^n$ is extremal (\cite{GHS:2}, \cite{GS:1}), we may apply this result to those ideals. It is an immediate consequence of the Eliahou--Kervaire theorem that if $I$ is a lex-segment ideal, a generic initial ideal, or the ideal of a $k$-configuration in $\P^n$ which has {\em no} generators in degree $d$, then $\beta_{q,i}=0$ whenever $i-q=d$.

\begin{remk}  \label{R:205} Let $I$ be any homogeneous ideal of $R=k[x_1,\dots,x_n]$ and $J=\Gin(I)$. Then, by Theorem~\ref{Cancellation Principle}, we have
$$
\beta_{q,i}(I) \le \beta_{q,i}(J).
$$
In particular, if $\beta_{q,i}(J)=0$, then $\beta_{q,i}(I)=0$.
\end{remk}

Let $I$ be a homogeneous ideal of $R=k[x_1,\dots,x_n]$ such that $\dim(R/I)=d$.
In \cite{HT}, they defined the {\em $s$-reduction number} $r_s(R/I)$ of $R/I$ for $s\ge d$ and have shown the following theorem.

\begin{thm}[\cite{AM}, \cite{HT}]\label{Reduction Number_2}
For a homogeneous ideal $I$ of $R$,
\[r_s(R/I)=r_s(R/\Gin(I)).\]
\end{thm}

If $I$ is a Borel fixed monomial ideal of $R=k[x_1,\dots,x_n]$
with $\dim(R/I)=n-d$, then we know that there are positive numbers
$a_1,\dots,a_d$ such that $x_i^{a_i}$ is a minimal generator of
$I$. In \cite{HT}, they have also proved that if a monomial ideal
$I$ is strongly stable, then
$$
r_s(R/I)=\min\{\ell \mid x^{\ell+1}_{n-s}\in I\}.
$$

Furthermore, the following useful lemma has been proved in
\cite{AM}.

\begin{lem}[Lemma 2.15, \cite{AM}] \label{Reduction Number 1}
For a homogeneous ideal $I$ of $R$ and for $s\geq \dim(R/I)$, the
$s$-reduction number $r_s(R/I)$ can be given as the following:
\begin{align*}
r_s(R/I)  = &\min\{\ell \,:\,x_{n-s}^{\ell+1}\in \Gin(I)\}\\
          = &\min\{\ell \,: \textup{
Hilbert function of }R/(I+J) \textup{ vanishes in degree }\ell+1
\textup{ }\}
\end{align*}
where $J$ is generated by $s$ general linear forms of $R$.
\end{lem}

For a homogeneous ideal $I$ of $R=k[x_1,\dots,x_n]$, we recall that $I^{\rm lex}$ is a lex-segment ideal associated to $I$. In Section~\ref{Sec:004}, we shall use the following two useful lemmas.

\begin{lem}\label{L:208}
Let $I$ be a homogeneous ideal of $R=k[x_1,\ldots,x_n]$ and let
$\bar{I}=(I_{\leq d+1})$ for some $d>0$. Then,

\begin{itemize}
\item[(a)] $\beta_{i,j}(I)\leq \beta_{i,j}(\Gin(I))\leq
\beta_{i,j}(I^{\rm lex})$ for all $i$, $j$.

\item[(b)]
$\beta_{0,d+2}(\bar{I}^{\rm lex})=\beta_{0,d+2}(I^{\rm lex})-\beta_{0,d+2}(I)$,

\item[(c)]
$\beta_{0,d+2}(\Gin(\bar{I}))=\beta_{0,d+2}(\Gin(I))-\beta_{0,d+2}(I).$
\end{itemize}
\end{lem}

\begin{proof}
(a) The first inequality can be proved by Theorem~\ref{Cancellation Principle}. The second one is directly obtained from
the theorem of  Bigatti, Hulett, and Pardue (\cite{Bi}, \cite{Hul}, and \cite{Pa}).

\medskip

(b) First note that
\begin{equation} \label{EQ:201}
\begin{array}{llllllllllllllllll}
&   & \beta_{0,d+2}(I^{\rm lex}) \\[1ex]
& = & \dim_k(I^{\rm lex})_{d+2}-\dim_k(R_1(I^{\rm lex})_{d+1})\\[1ex]
& = & [\dim_k R_{d+2}-\dim_k(R_1(I^{\rm lex})_{d+1})]
     -[\dim_k R_{d+2}-\dim_k(I^{\rm lex})_{d+2}]\\[1ex]
& = & \H_{R/I^{\rm lex}}(d+1)^{\langle d+1\rangle}-\H_{R/I^{\rm lex}}(d+2) \\[1ex]
& = & \H_{R/I}(d+1)^{\langle d+1\rangle}-\H_{R/I}(d+2) \qquad
(\because \H_{R/I}(t)=\H_{R/I^{\rm lex}}(t)  \text{ for every } t).
\end{array}
\end{equation}

It follows from equation~\eqref{EQ:201} that
$$
\begin{array}{llllllllllllllllllll}
 \beta_{0,d+2}(I)
& = & \dim_k(I_{d+2})-\dim_k(\bar{I}_{d+2})\\[1ex]
& = & [\dim_k R_{d+2} -\dim_k(\bar{I}_{d+2})]
-[\dim_k R_{d+2}- \dim_k(I_{d+2})] \\[1ex]
& = & \H_{R/\bar{I}}(d+2)-\H_{R/I}(d+2)\\[1ex]
& = & (\H_{R/I}(d+1)^{\langle d+1\rangle}-\H_{R/I}(d+2))
-(\H_{R/{I}}(d+1)^{\langle d+1\rangle}-\H_{R/\bar{I}}(d+2))\\[1ex]
& = & (\H_{R/I}(d+1)^{\langle d+1\rangle}-\H_{R/I}(d+2))-(\H_{R/{\bar I}}(d+1)^{\langle d+1\rangle}-\H_{R/\bar{I}}(d+2))\\[1ex]
&   & (\because \H_{R/I}(d+1)=\H_{R/{\bar I}}(d+1))\\[1ex]
& = & \beta_{0,d+2}(I^{\rm lex})-\beta_{0,d+2}(\bar{I}^{\rm lex})\qquad
(\because \text{equation}~\eqref{EQ:201}).
\end{array}
$$

\medskip

(c) Note that $\Gin(I)_{d+1}=\Gin(\bar{I})_{d+1}$. Hence we have
$$
\begin{array}{lllllllllllllllllll}
&   & \beta_{0,d+2}(I) \\[1ex]
& = & \dim_k(I_{d+2})-\dim_k(\bar{I}_{d+2})\\[1ex]
& = & \dim_k(\Gin(I)_{d+2})-\dim_k(\Gin(\bar{I})_{d+2})\\[1ex]
& = & [\dim_k(\Gin(I)_{d+2})-\dim_k(R_1\Gin(I)_{d+1})]-[\dim_k(\Gin(\bar{I})_{d+2})-\dim_k(R_1\Gin(\bar{I})_{d+1})] \\[1ex]
&   & (\because \Gin(I)_{d+1}=\Gin(\bar{I})_{d+1}) \\[1ex]
& = & \beta_{0,d+2}(\Gin(I))-\beta_{0,d+2}(\Gin(\bar{I})),
\end{array}
$$
which completes the proof.
\end{proof}

\begin{lem}\label{L:209}
Let $I \subset R=k[x_1,x_2x_3]$ be a homogenous ideal and let
that $A=R/I$ be a graded Artinian algebra. Then, for every $d>0$,
\begin{itemize}
\item[(a)]
$\beta_{1,d}(I^{\rm lex})-\beta_{1,d}(I)=[\beta_{0,d}(I^{\rm lex})-\beta_{0,d}(I)]+[\beta_{2,d}(I^{\rm lex})-\beta_{2,d}(I)].$

\item[(b)]
$\beta_{1,d}(\Gin(I))-\beta_{1,d}(I)=[\beta_{0,d}(\Gin(I))-\beta_{0,d}(I)]+[\beta_{2,d}(\Gin(I))-\beta_{2,d}(I)].$
\end{itemize}
\end{lem}

\begin{proof} It is immediate by the cancellation principle.
\end{proof}

\section{An $h$-vector of A Graded Artinian Level Algebra Having The WLP}

In this section, we think of $h$-vectors of a graded
Artinian level algebra with the WLP and we shall prove that some of
graded Artinian O-sequences are not level using generic initial
ideals. Moreover, we assume that $R=k[x_1,\dots,x_n]$ is an $n$-variable polynomial ring over a field $k$ with characteristic $0$.

\medskip

For positive integers $h$ and $i$, $h$ can be written uniquely in the form
$$
h=h_{(i)}:=\binom{m_i}{i}+\binom{m_{i-1}}{i-1}+\cdots+\binom{m_j}{j}
$$
where $m_i>m_{i-1}>\cdots >m_j\ge j\ge 1$. This expansion for $h$ is called the $i$-{\em binomial expansion} of $h$. For such $h$ and $i$, we define
$$
\begin{array}{llllllllllllllll}
(h_{(i)})^{-}
& := & \ds \binom{m_i-1}{i}+\binom{m_{i-1}-1}{i-1}+\cdots+\binom{m_j-1}{j}, \\[2.5ex]
(h_{(i)})^{+}_{+}
& := & \ds \binom{m_i+1}{i+1}+\binom{m_{i-1}+1}{i}+\cdots+\binom{m_j+1}{j+1}.
\end{array}
$$
Let $\H=\{h_i\}_{i\ge 0}$ be the Hilbert function of a graded ring $A$. For simplicity in the notation we usually rewrite
$((h_i)_{(i)})^{-}$ and $((h_i)_{(i)})^{+}_{+}$ as $(h_{i})^{-}$ and $(h_{i})^{+}_{+}$, respectively. Recall that we sometimes use another simpler notation $h^{\langle i\rangle}$ for  $(h_i)^{+}_{+}$ and define $0^{\langle i\rangle}=0$.

\medskip

A well known result of Macaulay is the following theorem.

\begin{thm}[Macaulay] Let $\H=\{h_i\}_{i\ge 0}$ be a sequence of non-negative integers such that $h_0=1$, $h_1=n$, and $h_i=0$ for every $i>e$. Then $\H$ is the $h$-vector of some standard graded Artinian algebra if and only if, for every $1\le d\le e-1$,
$$
h_{d+1}\le (h_d)^{+}_{+}=h_d^{\langle d\rangle}.
$$
\end{thm}

We use a generic initial ideal with respect to the reverse lexicographic
order to obtain results in Section 3. Note that, by Green's hyperplane restriction theorem (see \cite{ERV}), we have that
\begin{equation} \label{EQ:310}
\H({R/(J+x_n)},d)\le (\H({R/J},d))^{-},
\end{equation}
and the equality holds when $J$ is a strongly stable ideal of $R$. In particular, the equality holds for any lex-segment ideal since a lex-segment ideal $J$ of $R=k[x_1,\dots,x_n]$ is  also a strongly stable ideal.

\medskip

The following lemma will be used often in this section.

\begin{lem}\label{lemma 1}
Let $A=R/I$ be an Artinian $k$-algebra and let $L$ be a general
linear form.
 \begin{enumerate}
  \item[(a)] If
   \[\dim_k(0:L)_d > (n-1)\dim_k(0:L)_{d+1}\]
    for some $d>0$, then $A$ has a socle element in degree $d$.
  \item[(b)] Let $h(A)=(h_0,h_1,\dots,h_s)$ be the {\em$h$-vector} of
  $A$. Then, we have
  \begin{equation}\label{eq201}
  h_d-h_{d+1}\leq \dim_k(0:L)_d \leq h_d-h_{d+1}+(h_{d+1})^{-}.
  \end{equation}
   In particular, $\dim_k(0:L)_d =h_d-h_{d+1}$ if and only if $d\geq r_1(A)$.
 \end{enumerate}
\end{lem}

\begin{proof} (a) Consider a map $\varphi: (0:L)_d \rightarrow \bigoplus^{n-1}
 (0:L)_{d+1}$, defined by $\varphi(F)=(x_1F,\ldots,x_{n-1}F)$.
 Since $L$ is a general linear form, we may assume that the kernel of this map is
 exactly $\soc(A)_{d}$. Since $\dim_k(0:L)_d > (n-1)\dim_k(0:L)_{d+1}$, the map
 $\varphi$ is not injective and we obtain the desired result.

\medskip

 (b) Consider the following exact sequence
   \[0 \rightarrow (0:L)_d
    \rightarrow A_{d} \overset{\times L}{\rightarrow} A_{d+1} \rightarrow
    (A/LA)_{d+1}  \rightarrow 0.\]
 Then we have
    \begin{equation}\label{eq202}
    \dim_k(0:L)_d= h_d-h_{d+1}+\dim_k [A/(L)A]_{d+1},
    \end{equation}
 and thus $h_d-h_{d+1}\leq \dim_k(0:L)_d $. The right hand side of
 the inequality (\ref{eq201}) follows from Green's hyperplane restriction theorem,
 i.e., $\dim_k [A/(L)A]_{d+1}\leq (h_{d+1})^{-}$.

\medskip

Moreover, $\dim_k(0:L)_d =h_d-h_{d+1}$ if and only if $\dim_k [A/(L)A]_{d+1}=0$,
     and it is equivalent to $d\geq r_1(A)$ by the definition of $r_1(A)$.
\end{proof}

\begin{remk}\label{R:302}
Let $\H=(h_0,h_1,\ldots,h_s)$ be the $h$-vector of a graded
Artinian level algebra $A=R/I$ and $L$ is a general linear form of $A$. In general, it is not easy to find the reduction number $r_1(A)$ based on its $h$-vector.
However
if $h_{d+1}\leq d+1$ then $(h_{d+1})^{-}=0$, and thus $\dim_k(0:L)_d =h_d-h_{d+1}$. Hence $d\geq r_1(A)$ by Lemma~\ref{lemma 1}. In other words,
$$
r_1(A)\leq \min\{\,k\,\mid h_{k+1}\,\leq {k+1}\}.
$$
\end{remk}

\begin{prop}\label{Proposition 0}
Let $R=k[x_1,\ldots,x_n]$ and let $\H=(h_0,h_1,\ldots,h_s)$ be the
$h$-vector of a graded Artinian level algebra $A=R/I$ with socle
degree $s$. Suppose that $h_{d-1}>h_{d}$ for some $d\geq r_1(A)$.
Then
 \begin{enumerate}
  \item[(a)] $h_{d-1}>h_{d}>\cdots>h_{s-1}>h_s>0$, and
  \item[(b)] $h_{t-1}-h_t \leq (n-1)(h_t-h_{t+1})$ for all $d\le t\le s$.
 \end{enumerate}
\end{prop}

\begin{proof}
(a) First of all, note that, by Lemma~\ref{lemma 1} (b), $h_{t}-h_{t+1}=\dim_k(0:L)_{t}$ for every $t\ge r_1(A)$. Hence we have that
$$
h_{d-1}>h_d\ge h_{d+1}\ge \cdots \ge h_s.
$$
Now assume that there is $t\geq d$ such
that $h_{t-1}>h_{t}=h_{t+1}$. Since $t\geq r_1(A)$, we know that,
by Lemma~\ref{lemma 1} (b),
$$
\dim_k(0:L)_{t-1}\ge h_{t-1}-h_t> 0 \quad \textup{ and } \quad \dim_k(0:L)_{t}=0.
$$
Hence there is a socle element of $A$ in degree $t-1$, which is a contradiction since $A$ is level. This means that $h_t>h_{t+1}$ for every $t\ge d-1$.

\medskip

(b) Since $A$ is a level algebra and
$\dim_k(0:L)_{t}=h_{t-1}-h_t$, the result follows directly from
Lemma~\ref{lemma 1} (a).
\end{proof}

\begin{remk} \label{R:035}
Let $I$ be a homogeneous ideal of $R=k[x_1,\dots,x_n]$ such that $R/I$ has the WLP with a Lefschetz element $L$ and let $\H(R/I,d-1)>\H(R/I,d)$ for some $d$.
Now we consider the following exact sequence
\begin{equation} \label{EQ:306}
(R/I)_{d-1} \overset{\times L}{\rightarrow} (R/I)_d
\rightarrow (R/(I+(L)))_{d} \rightarrow 0.
\end{equation}
Since $R/I$ has the WLP and $\H(R/I,d-1)>\H(R/I,d)$, the
above multiplication map cannot be injective, but surjective. In
other words, $(R/(I+(L)))_{d} = 0$. This implies that
$d>r_1(R/I)$ by Lemma~\ref{Reduction Number 1}.
\end{remk}

The following theorem shows a useful condition to be a level O-sequence with the WLP.

\begin{thm} \label{T:307}
Let $R=\kxn{n}$, $n\geq 3$ and let $\H=(h_0,h_1,\dots,h_s)$ be the
Hilbert function of a graded Artinian level algebra $A=R/I$ having
the WLP. Then,
\begin{itemize}
\item[(a)] the Hilbert function $\H$ is a strictly unimodal
O-sequence
$$
h_0<h_1<\cdots<h_{r_1(A)}=\cdots=h_{\theta}>\cdots>h_{s-1}>h_s
$$
 such that the positive part of the first difference
$\Delta\H$ is an O-sequence, and

\medskip

\item[(b)] $
\begin{array}{rllllllllllllllllll}
 h_{d-1}-h_d \leq (n-1)(h_d-h_{d+1})
\end{array}
$ for $s \geq d>\theta$.
\end{itemize}
\end{thm}

\begin{proof}
(a) First, note that, by Proposition 3.5 in \cite{HMNW},
$\H$ is a unimodal O-sequence such that the positive part of the
first difference is an O-sequence. Hence it suffices to show that $\H$ is strictly unimodal.

If $d\le r_1(A)$, then $\H_{R/(I+L)}(d)\ne 0$ by the definition of $r_1(A)$, and so the multiplication map $\times L$ is not surjective in equation~\eqref{EQ:306}. In other words, the multiplication map $\times L$ is injective since $A$ has the  WLP. Thus, we have a short exact sequence as follows
$$
0 \rightarrow (R/I)_{d-1} \overset{\times L}{\rightarrow} (R/I)_d
\rightarrow (R/(I+(L)))_{d} \rightarrow 0.
$$
Hence we obtain that
$$
\begin{array}{lllllllllllllllllllllllllllll}
\H_A(d)
& = & \H_A(d-1)+\H_{R/(I+L)}(d) \\[1ex]
& > & \H_A(d-1) \qquad (\because \H_{R/(I+L)}(d)\ne 0),
\end{array}
$$
and so the Hilbert function of $A$ is strictly increasing up to $r_1(A)$.

Moreover, by Proposition \ref{Proposition 0} (a), $\H$ is strictly decreasing in degrees $d\geq \theta$, where
$$
\theta:=\min\{ t \mid h_t>h_{t+1}\}.
$$

(b) The result follows directly from Proposition \ref{Proposition
0} (b).
\end{proof}

\begin{remk}
Theorem \ref{T:307} gives us a necessary condition when a
numerical sequence becomes a level O-sequence with the WLP. In general,
this condition is not sufficient. One can find many non-level
sequences satisfying the inequality of Theorem \ref{T:307} in
\cite{GHMS}.
\end{remk}

In \cite{GHMS}, they gave some of ``non-level sequences" using
homological method, which is the combinatorial notion of
the cancellation of shifts in the minimal free resolutions of the lex-segment ideals associated to the given homogenous ideals.

\medskip

In this section, we shall use generic initial ideals, instead of the lex-segment ideals. First note that, by Bigatti-Hulett-Pardue Theorem,  the worst minimal free resolution of a homogenous ideal $I$ depends on only the Hilbert function of $I$.
Unfortunately, we cannot apply their theorem to obtain the minimal free resolutions of  the generic initial ideals. However, we can find Betti-numbers $\beta_{i,\,d+i}({\rm Gin}(I))$ for $d > r_1(A)$ and $i\ge 0$, which depends on only the given Hilbert function (see Corollary~\ref{C:303}).

\medskip

For the rest of this section, we need the following useful results.

\begin{lem}\label{strongly stable}
Let $J$ be a stable ideal of $R$ and let $T_1,\ldots,T_r$ be the
monomials which form a $k$-basis for
$\left((J:x_n)/J\right)_{d-1}$ then
$$
\{x_nT_1,\ldots,x_nT_r\}=\{T\in \mathcal G(J)_d \mid \,x_n\text{
divides } T\, \}.
$$
In particular,
$$
\dim_k \left((J:x_n)/J\right)_{d-1}= \big|\{T\in \mathcal G(J)_d
\mid \,x_n\text{ divides } T\, \}\big|.
$$
\end{lem}

\begin{proof}
For every $T=x_nT^{'}\in \mathcal G(J)_d$, we have that $x_nT'\in
J_d\subset J$, i.e., $T'\in (J:x_n)_{d-1}$, and thus
$\overline{T}'\in ((J:x_n)/J)_{d-1}=\langle
\overline{T}_1,\dots,\overline{T}_r\rangle$. However, since $T'$
and $T_i$ are all monomials of $(J:x_n)_{d-1}$ in degree $d-1$, we
have that $T'=T_i$ for some $i$, and hence $T=x_nT'\in
\{x_nT_1,\ldots,x_nT_r\}$.

Conversely,  note that $T_i\notin J_{d-1}$ and $x_nT_i\in J_d$ for
every $i=1,\dots,r$. If $x_nT_i\notin {\mathcal G}(J)_d$ for some
$i=1,\dots,r$, then $x_nT_i\in R_1J_{d-1}$. Since $T_i\notin
J_{d-1}$, we see that
$$
x_nT_i=x_jU
$$
for some monomial $U\in J_{d-1}$ and $j<n$. Hence we have that
$$
x_n\mid U.
$$
Moreover, since $J$ is a stable monomial ideal, for every $\ell
<n$,
$$
\frac{x_\ell}{x_n}U\in J_{d-1}.
$$
In particular, we have
$$
T_i=\frac{x_j}{x_n}U\in J_{d-1},
$$
which is a contradiction. Therefore,  $x_nT_i\in {\mathcal
G}(J)_d$, for every $i=1,\dots,r$, as we wished.
\end{proof}

Using the previous lemma, we obtain the following proposition,
which we know the difference between $h_d$ and $h_{d+1}$ when
$d>r_1(A)$.

\begin{prop}\label{Proposition 1}
Let $A=R/I$ be a graded Artinian algebra with Hilbert function
$\bh=(h_0,h_1,\ldots,\linebreak h_s)$ and let $J=\Gin(I)$. If
$d\geq r_1(A)$ then,
$$\big|\{T\in \mathcal G(J)_{d+1} \mid \,x_n\text{
divides } T\, \}\big|=h_d-h_{d+1}.$$ Moreover, if $d>r_1(A)$,
$$
|\,\mathcal G(J)_{d+1}|=\big|\{T\in \mathcal G(J)_{d+1} \mid
\,x_n\text{ divides } T\, \}\big|=h_d-h_{d+1}.
$$
\end{prop}

\begin{proof}
Consider the following exact sequence:
$$
0 \rightarrow \left ((J:x_n)/J \right )_{d} \rightarrow (R/J)_{d}
\stackrel{\times x_n}{\longrightarrow} (R/J)_{d+1} \rightarrow
(R/J + (x_n))_{d+1} \rightarrow 0.
$$
Note that $\H(R/I,t)=\H(R/J,t)$ for every $t\ge 0$. So
\begin{equation} \label{EQ:031}
\begin{array}{rclcccccccccccccc}
\dim_k \left ((J:x_n)/J \right )_{d}+\dim_k (R/J)_{d+1}
& = & \dim_k (R/J)_{d}+ \dim_k (R/J + (x_n))_{d+1}, \\[1ex]
\Leftrightarrow\quad \dim_k \left ((J:x_n)/J \right )_{d}+h_{d+1}
& = & h_{d}+ \dim_k (R/J + (x_n))_{d+1}.
\end{array}
\end{equation}
Moreover, by Theorem~\ref{Generic initial ideal of hyperplane
section}, Theorem~\ref{Reduction Number_2}, and
Lemma~\ref{Reduction Number 1}, we have
$$
\begin{array}{lllllllllllllllllll}
r_1(R/I)
& = & r_1(R/J) \\
& = &   \min\{\,\ell \mid  \H(R/J+(x_n),\ell+1)=0\, \},
\end{array}
$$
which means $\H(R/J+(x_n),d+1)=0$ for every $d\geq r_1(R/I)$.
Hence, from equation~\eqref{EQ:031}, we obtain
\begin{equation} \label{EQ:032}
\dim_k \left ((J:x_n)/J \right )_{d} = \big|\{T\in \mathcal
G(J)_{d+1} \mid \,x_n\text{ divides } T\, \}\big| = h_{d}-h_{d+1}.
\end{equation}
Now suppose that $d>r_1(A)$. Then it is obvious that
\begin{equation} \label{EQ:033}
\{T\in \mathcal G(J)_{d+1}\mid \,x_n\text{ divides } T\, \}
\subseteq \mathcal G(J)_{d+1}.
\end{equation}

Conversely, note that $x_{n-1}^{d}\in J$ from the first equality
of Lemma~\ref{Reduction Number 1}. Since $J$ is a strongly stable
ideal, $J_d$ has to contain all monomials $U$ of degree $d$ such
that
$$
{\rm supp}({U}):=\{i \mid x_i \text{ divides } U\} \subseteq
\{1,\dots,n-1\}.
$$
This implies $\overline {\bf m}_d\subseteq J_d$ where $\overline
{\bf m}=(x_1,\ldots,x_{n-1})^d$. Thus we have
$$
R_1 \overline {\bf m}_d\subseteq J_{d+1}.
$$
Therefore, for every $T\in \mathcal G(J)_{d+1}$, we have $x_n \mid
T$, and so
\begin{equation} \label{EQ:034}
\mathcal G(J)_{d+1} \subseteq \{T\in \mathcal G(J)_{d+1}\mid
\,x_n\text{ divides } T\, \}.
\end{equation}
It follows from equations~\eqref{EQ:033} and \eqref{EQ:034} that
\begin{equation} \label{EQ:035}
\mathcal G(J)_{d+1} = \{T\in \mathcal G(J)_{d+1}\mid \,x_n\text{
divides } T\,\},
\end{equation}
and hence
$$
\mid\mathcal G(J)_{d+1}\mid=\dim_k \left ((J:x_n)/J \right )_{d} =
h_{d}-h_{d+1},
$$
as we wished.
\end{proof}

\begin{cor}\label{C:303}
Let $A=R/I$ be a graded Artinian algebra with Hilbert function
$\bh=(h_0,h_1,\ldots,\linebreak h_s)$.  If
$d>r_1(A)$ then, for all $i\geq 0$,
\[\beta_{i,\,i+(d+1)}(\Gin(I))=(h_{d}-h_{d+1})\binom{n-1}{i}. \]
\end{cor}

\begin{proof}  By Proposition \ref{Proposition 1},
$$|\,\mathcal
G(\Gin(I))_{d+1}|=\big|\{T\in \mathcal G(\Gin(I))_{d+1} \mid \,x_n\text{
divides } T\, \}\big|=h_{d}-h_{d+1}
$$
for every $d>r_1(A)$, and thus the result follows from Theorem \ref{EK}.
\end{proof}

Recall that a homogeneous ideal $I$ is $m$-{\em regular} if, in the minimal free resolution of $I$, for all $p \ge 0$, every $p$-th syzygy has degree $\le m + p$. The regularity of $I$, ${\rm reg}(I)$, is the smallest such $m$.

In \cite{BS} and \cite{G}, they proved that the regularity of $\Gin(I)$ is the largest degree of a generator of $\Gin(I)$. Moreover, Bayer and Stillman \cite{BS} showed the regularity of $I$ is equal to the regularity of $\Gin(I)$.

\begin{thm}[\cite{BS}, \cite{G}] \label{T:311-2}
For any homogeneous ideal $I$, using the reverse
lexicographic order,
$$
{\rm reg}(I) = {\rm reg}(\Gin(I)).
$$
\end{thm}

\begin{thm}[Crystallization Principle, \cite{AM} and \cite{G}] \label{CP}
Let $I$ be a homogeneous ideal generated
in degrees $\le d$. Assume that there is a monomial order $\tau$ such that $\Gin_\tau(I)$ has no generator in degree $d + 1$. Then $\Gin_\tau(I)$ is generated in degrees $\le d$ and $I$ is $d$-regular.
\end{thm}

\begin{lem}\label{lem}
Let $R=k[x_1,x_2,x_3]$ and let $A=R/I$ be an Artinian algebra and
let $\H=(h_0,h_1,\dots,h_s)$ be the Hilbert function of $A=R/I$.
Suppose that, for $t>0$,
\begin{enumerate}
\item[(a)] $\soc (A)_{t-2}=0$,
\item[(b)]
$\beta_{1,\,t+1}(\Gin(I))=\beta_{2,\,t+1}(\Gin(I))$.
\end{enumerate}
Then $(I_{\leq t})$ is $t$-regular and
\begin{equation}\label{eq:307}
h_{t-1}-h_{t}\leq \dim_k\soc (A)_{t-1}\leq
h_{t-1}-h_{t}+(h_{t})^{-}
\end{equation}

In particular, if $t> r_1(A)$ then
$$
\dim_k(\soc (A)_{t-1}) = h_{t-1}-h_{t}.
$$
\end{lem}

\begin{proof}
Let $\bar I=(I_{\leq t})$. Note that $\beta_{i,t+1}(\Gin(I))=\beta_{i,t+1}(\Gin(\bar I))$ for $i=1,\, 2$ and $\beta_{0,t+1}(\bar I)=0$. Furthermore, since $I$ and $\bar I$ agree in degree $\le t$ and $\soc (A)_{t-2}=0$, we see that $\beta_{2,t+1}(I)=\beta_{2,t+1}(\bar I)=0$.

Applying Lemma
\ref{L:209} (b)  the ideal $\bar I$, we have that
$$
\begin{array}{lllllllllllllllllllll}
& \beta_{1,t+1}(\Gin(\bar I))-\beta_{1,t+1}(\bar I)=
(\beta_{0,t+1}(\Gin(\bar I))-\beta_{0,t+1}(\bar I))+
(\beta_{2,t+1}(\Gin(\bar I))-\beta_{2,t+1}(\bar I)) \\
\Rightarrow &
-\beta_{1,t+1}(\bar I)=
(\beta_{0,t+1}(\Gin(\bar I))-\beta_{0,t+1}(\bar I))
-\beta_{2,t+1}(\bar I) \quad (\because \beta_{1,\,t+1}(\Gin(\bar I))=\beta_{2,\,t+1}(\Gin(\bar I))) \\
\Rightarrow &
-\beta_{1,t+1}(\bar I)=\beta_{0,t+1}(\Gin(\bar I))  \quad
(\because
\beta_{0,t+1}(\bar I)=\beta_{2,t+1}(\bar I)=0) \\
\Rightarrow &
\beta_{0,\,t+1}(\Gin(\bar I))=0.
\end{array}
$$
Thus, by Theorem \ref{CP}, the ideal $\bar I=(I_{\leq t})$ is $t$-regular.

\medskip

Let $\bar A=R/\bar I$. For a general linear form $L$, consider the following exact
sequence
\begin{equation} \label{EQ:311-1}
0 \rightarrow \left (0:_{\bar A} L \right )_{t-1} \rightarrow (R/\bar
I)_{t-1} \stackrel{\times L}{\longrightarrow} (R/\bar I)_{t}
\rightarrow (R/\bar I + (L))_{t} \rightarrow 0.
\end{equation}
After we replace $\bar I$ and $\bar A$ by $\Gin(\bar I)$ and $\tilde  A=R/\Gin(\bar I)$, respectively, we can rewrite equation~\eqref{EQ:311-1} as
\begin{equation} \label{EQ:312-1}
0 \rightarrow \left (0:_{\tilde A} x_3 \right )_{t-1} \rightarrow (R/\Gin(\bar
I))_{t-1} \stackrel{\times x_3}{\longrightarrow} (R/\Gin(\bar I))_{t}
\rightarrow (R/\Gin(\bar I) + (x_3))_{t} \rightarrow 0.
\end{equation}
Then, by Theorem \ref{Generic initial ideal of
hyperplane section}, we know that

\begin{align*}
\dim_k \left  (0:_{\tilde A} x_3 \right )_{t-1}
  &= \dim_k((\Gin(\bar I):x_3)/\Gin(\bar I))_{t-1} \\
  &= h_{t-1}-h_{t}+\dim_k (R/\Gin(\bar I) + (x_3))_{t}\\
  &= h_{t-1}-h_{t}+\dim_k \left (R/\bar I + (L)\right )_{t}\\
  &= \dim_k (0:_{\bar A} L)_{t-1}.
\end{align*}
On the other hand, by Lemma~\ref{strongly stable},
$$
\begin{array}{llllllllllllllllll}
\dim_k((\Gin(\bar I):x_3)/\Gin(\bar I))_{t-1}
& = & \big|\{ T\in \GG(\Gin(\bar I))_{t} \mid x_3 \text{ divides } T\}\big| \\
& = & \beta_{2,\,t+2}(\Gin(\bar I)),
\end{array}
$$
and by Lemma~\ref{lemma 1} (b)
\begin{equation}\label{EQ:313-1}
h_{t-1}-h_{t}\leq \dim_k ((0:_{\bar A} L)_{t-1})\leq
h_{t-1}-h_{t}+(h_{t})^{-}.
\end{equation}
Note that, by Theorem~\ref{CP}, $\beta_{1,t+2}(\Gin(\bar I))=0$ since $\bar I=(I_{\leq t})$ is $t$-regular. Moreover, since $I$ and $\bar I$ agree in degree $\le t$, we have $\beta_{2,t+2}(I)=\beta_{2,t+2}(\bar I)$. Hence, by Theorem
\ref{Cancellation Principle},
\begin{equation} \label{EQ:314-1}
  \begin{array}{llllllllllllllllll}
  \dim_k\soc(A)_{t-1}
  &=\beta_{2,\,t+2}(I)\\
  &=\beta_{2,\,t+2}(\bar I)\\
  &=\beta_{2,\,t+2}(\Gin(\bar I)) \quad (\because \beta_{1,t+2}(\Gin(\bar I))=0)\\
  &= \dim_k (0:_{\bar A} L)_{t-1}.
\end{array}
\end{equation}
Hence it follows from equations~\eqref{EQ:313-1} and ~\eqref{EQ:314-1}, we obtain the inequality~\eqref{eq:307}. Moreover, by Lemma~\ref{lemma 1} (b), we have
$$
\dim_k(\soc (A)_{t-1}) = h_{t-1}-h_{t} \quad \textup{ for } \quad t>r_1(A),
$$
as we wished.
\end{proof}

\begin{thm}\label{T:301-1}
Let $A=R/I$ be an Artinian algebra of codimension $3$ with socle
degree $s$. If
 \begin{equation}\label{eq302}
   \beta_{1,\,d+2}(\Gin(I))=\beta_{2,\,d+2}(\Gin(I))>0.
 \end{equation}
for some $d<s$, then $A$ is not level.
\end{thm}

\begin{proof} Assume $A$ is level. Then $\beta_{2,d+2}(I)=\soc(A)_{d-1}=0$, and hence, by Lemma~\ref{lem}, $\bar I=(I_{\leq d+1})$ is $(d+1)$-regular.

\medskip

Let $\bar A=R/\bar I$. Note that $\soc(A)_{d}=\soc(\bar A)_{d}$ since $A$ and $\bar A$ agree in degree $\le d+1$, i.e.,
$$
\dim_k\soc(A)_d=\beta_{2,d+3}(I)=\beta_{2,d+3}(\bar I)=\dim_k\soc(\bar A)_d.
$$

For a general linear form $L$, by Lemmas~\ref{lemma 1}~(a)
and~\ref{strongly stable}, we have that
$$
\begin{array}{llllllllllllllllll}
0
& < & \beta_{2,\,d+2}(\Gin(I)) & (\because \text{ by assumption})\\[1ex]
& = & \ds \sum_{T\in \GG(\Gin(I))_d}\binom{m(T)-1}{2} \\[4ex]
& = & \dim_k \left[(\Gin(I):x_3)/\Gin(I)\right]_{d-1} & (\because \text{ by Lemma~\ref{strongly stable}})\\[1ex]
& = & \dim_k\left[(I:L)/I)\right]_{d-1}\\[1ex]
& \leq & 2\dim_k\left[(I:L)/I)\right]_{d} & (\because \text{ by Lemma~\ref{lemma 1} (a) and $\soc(A)_{d-1}=0$}).
\end{array}
$$
Note that, in the similar way, we have
$\beta_{2,\,d+3}(\Gin(I))=\dim_k\left[(I:L)/I)\right]_{d}$. Hence
\[\beta_{2,\,d+3}(\Gin(I))>0.\]
Since $\bar I=(I_{\leq
d+1})$ is $(d+1)$-regular and ${\rm reg}(\bar I)={\rm
reg}(\Gin(\bar I))$ by Theorem~\ref{T:311-2}, we have that
$$
\begin{array}{lllllllllllllllllllllllllll}
            &\beta_{0,\,d+3}(\Gin(\bar I))=\beta_{1,\,d+3}(\Gin(\bar I))=0,\\
            & \beta_{0,\,d+3}(\bar I)= \beta_{1,\,d+3}(\bar I)=0.
\end{array}
$$
Thus, by Lemma \ref{L:209} (b),
$$
\beta_{2,\,d+3}(\bar I)=\beta_{2,\,d+3}(\Gin(\bar I))>0,
$$
which follows that $R/\bar I$ has a socle element in degree $d$, so does $R/I$. This is a contradiction, and thus we complete the proof.
\end{proof}

\begin{rem}
Now we shall show that there is a level O-sequence satisfying Theorem \ref{T:307} (a) and (b), but it cannot be the Hilbert function of an Artinian algebra with the WLP.

\medskip

Consider an $h$-vector $\H =(\begin{matrix} 1 , 3 , 6 , 10 , 8 , 7 \end{matrix})$,
which was given in \cite{GHMS}. Furthermore, it has been shown that there is a level algebra of codimension $3$ with Hilbert function $\H$ in \cite{GHMS}.
They also raised a question if there exists a codimension 3 graded level algebra having the WLP with  Hilbert function $\H$.
Note that this is a codimension $3$ level O-sequence which satisfies the
condition in Theorem \ref{T:307}.

\medskip

Now suppose that there is an Artinian level algebra $A=R/I$ having the WLP with Hilbert function $\H$. In \cite{GHMS}, they gave several results about level or non-level sequences of graded Artinian algebras. One of the tools they used
was the fact that Betti numbers of a homogeneous ideal $I$
can be obtained by cancellation of the Betti numbers of
$I^{\lex}$. However, in this case, it is not available if $\H$ can be the Hilbert function of  an Artinian level algebra having the WLP based on the Betti numbers of $I^{\rm lex}$.

In fact, the  Betti diagram of $R/I^{\lex}$  is

\begin{center}
  \begin{verbatim}
                        total: 1   -   -   -
                        ------------------------
                            0: 1   -   -   -
                            1: 0   0   0   0
                            2: 0   0   0   0
                            3: 0   7   9   3
                            4: 0   2   4   2
                                    ......
  \end{verbatim}
\end{center}
and thus we cannot decide if there is a socle element of $R/I$
in degree $3$.

\medskip

Note that, by Theorem~\ref{T:307}, $r_1(A)=3$ since $A$ has the WLP. Hence,  by Corollary~\ref{C:303},
$$
\begin{array}{lllllllllllllll}
\beta_{2,6}(\Gin(I))
& = & (h_4-h_5)\binom{2}{2}=2\cdot 1=2, \quad \text{and} \\[1ex]
\beta_{1,6}(\Gin(I))
& = & (h_5-h_6)\binom{2}{1}=1\cdot 2=2.
\end{array}
$$
Therefore, by Theorem~\ref{T:301-1}, there is a socle element in $A$ in degree $3$, which is a contradiction. In other words, any Artinian level algebra $A$ with Hilbert function $\H$ does not have the WLP.
\end{rem}

\begin{remk}\label{R:301}
In general, Theorem~\ref{T:301-1} is not true if equation
\eqref{eq302} holds in the socle degree. For example, we
consider a Gorenstein sequence
\begin{center}
  \begin{tabular}{l|cccccccccccccccccccccccccccccc}
   $d$   & 0 & 1 & 2 & 3  & 4   \\
   \hline $h_d$ & 1 & 3 & 6 & 3 & 1
  \end{tabular}
 \end{center}

\ni
By Remark~\ref{R:302}, $r_1(A)\le 2$. Hence
$$
\begin{array}{llllllllllllllll}
\beta_{1,6}(\Gin(I))=(h_4-h_5)\binom{2}{1}=1\cdot 2=2, \quad \text{and}\\[1ex]
\beta_{2,6}(\Gin(I))=(h_3-h_4)\binom{2}{2}=2\cdot 1=2.
\end{array}
$$
Note this satisfies the condition of Theorem~\ref{T:301-1} in the socle degree,
but it is a level sequence.
\end{remk}

\begin{rem}
Let $A=R/I$ be an Artinian algebra and let
$\H=(h_0,h_1,\dots,h_s)$ be the Hilbert function of $A=R/I$. Then
an ideal $(I_{\leq d+1})$ is $(d+1)$-regular, if the Hilbert
function $\H$ of $A$ has the maximal growth in degree $d>0$, i.e.
$h_{d+1}=h_d^{\langle d\rangle}$. In particular, if $h_d = h_{d+1} = \ell \leq d,
$
then we know that $(I_{\leq d+1})$ is $(d+1)$-regular. Recently,
this result was improved in \cite{AM}, that is, $(I_{\leq d+1})$
is $(d+1)$-regular if
$
h_d = h_{d+1} \textup{ and } r_1(A)< d.
$

Note that, by Lemma~\ref{lemma 1}, the $k$-vector space dimension of $(0:L)_d$ in
degree $d\geq r_1(A)$ is $h_d-h_{d+1}$. By Proposition
\ref{Proposition 0}, we have a bound for the growth of Hilbert
function of $(0:L)$ in degree $d\geq r_1(A)$ if an Artinian
algebra $A$ has no socle elements in degree $d$. Theorem
\ref{T:302} shows that a similar result still holds on the
maximal growth of the Hilbert function of $(0:L)$ in codimension
three case.
\end{rem}

\begin{lem} \label{L:316-1}
Let $R=k[x_1,\ldots,x_n]$ and
let $A=R/I$ be an Artinian algebra with an $h$-vector $\H=(1,3,h_2,\cdots,h_s)$. If
$
h_{d-1}-h_d = (n-1)(h_d-h_{d+1})
$
for $r_1(A)< d < s$, then
$$
\beta_{(n-1),(n-1)+d}(\Gin(I))=\beta_{(n-2),(n-1)+d}(\Gin(I)).
$$
\end{lem}

\begin{proof} Let $J=\Gin(I)$. By Proposition~\ref{Proposition 1}, we have that
\begin{align*}
\beta_{(n-1),(n-1)+d}(J) & = \sum_{T\in \, \mathcal G(J)_d}\binom{m(T)-1}{n-1} \\
                     & = h_{d-1}-h_d.
\end{align*}
Moreover, by Corollary~\ref{C:303},
\begin{align*}
\beta_{(n-2),(n-1)+d}(J)
& =\beta_{(n-2),(n-2)+(d+1)}(J) \\
& = (h_d-h_{d+1})\binom{n-1}{n-2} \\
 &=(n-1)(h_d-h_{d+1})\\
& = h_{d-1}-h_d \qquad \qquad (\because \text{ by given condition})\\
& = \beta_{(n-1),(n-1)+d}(J),
\end{align*}
as we wished.
\end{proof}

\begin{thm}\label{T:302}
Let $R=k[x_1,x_2,x_3]$ and let $A=R/I$ be an Artinian algebra with an $h$-vector $\H=(1,3,h_2,\cdots,h_s)$.
If $\soc(A)_{d-1}=0$ and the Hilbert function of $(0:L)$ has a maximal growth in degree $d$ for $r_1(A) <d<s$, i.e.,
$
h_{d-1}-h_d = 2(h_d-h_{d+1}),
$
then
\begin{enumerate}
\item[(a)] $(I_{\leq\,d+1})$ is $(d+1)$-regular, and
\item[(b)] $\dim_k \soc(A)_d=h_d-h_{d+1}$.
\end{enumerate}
\end{thm}

\begin{proof} By Lemma~\ref{L:316-1}, we have
 \begin{equation}\label{eq301-2}
   \beta_{1,\,d+2}(\Gin(I))=\beta_{2,\,d+2}(\Gin(I)),
 \end{equation}
for $r_1(A) <d<s$, and  the result immediately follows from Lemma~\ref{lem}.
\end{proof}

\begin{cor}\label{C:302}
Let $R=k[x_1,x_2,x_3]$ and let $A=R/I$ be an Artinian algebra with an $h$-vector $\H=(1,3,h_2,\cdots,h_s)$. If
$h_{d-1}-h_d = 2(h_d-h_{d+1})>0$ for $r_1(A) <  d <s$, then $A$ is not level.
\end{cor}

\begin{proof}  By Lemma~\ref{L:316-1}, we have
$$
\beta_{2,d+2}(\Gin(I))=\beta_{1,d+2}(\Gin(I))>0,
$$
and hence, by Theorem~\ref{T:301-1}, $A$ cannot be level, as we wanted.
\end{proof}

\begin{remk}
Remark \ref{R:301} shows Corollary \ref{C:302} is not true if
$d=s$. However we know $h_{s-1}\leq 3h_s$ by Theorem \ref{T:307}.
\end{remk}

\begin{exmp}
Let $A=R/I$ be a codimension $3$ Artinian algebra and let $r_1(A)< d < s$. If $A$ has the Hilbert function
 \begin{center}
  \begin{tabular}{l|cccccccccccccccccccccccccccccc}
   $d$   & $\cdots$ & $d-1$ & $d$ & $d+1$  & $\cdots $  \\
   \hline $h_d$ & $\cdots$ & $a+3k$ & $a+k$ & $a$ & $\cdots$
  \end{tabular}
 \end{center}
such that $a>0$ and $k>0$, then by Corollary~\ref{C:302} $A$ cannot be level since
$$
h_{d-1}-h_d=2k=2(h_d-h_{d+1})
\Leftrightarrow \beta_{2,d+2}(\Gin(I))=\beta_{1,d+2}(\Gin(I))>0.
$$
\end{exmp}

For the codimension $3$ case, we have the following theorem, which follows from Theorems~\ref{T:307} and \ref{T:302} and
Corollary \ref{C:302}, and so we shall omit the proof here.

\begin{thm}\label{T:303}
Let $A=R/I$ be a graded Artinian level algebra of codimension 3
with the WLP and let $\H=(h_0,h_1,\ldots,h_s)$ be the Hilbert function
of $A$. Then,
 \begin{enumerate}
  \item[(a)] the Hilbert function $\H$ is a strictly unimodal
  $O$-sequence
  $$
h_0<h_1<\cdots<h_{r_1(A)}=\cdots=h_{\theta}>\cdots>h_{s-1}>h_s
$$
 such that the positive part of the first difference
$\Delta\H$ is an O-sequence, and
  \item[(b)] $h_{d-1}-h_d < 2(h_d-h_{d+1})$ for $s>d>\theta$.
  \item[(c)] $h_{s-1}\leq 3h_s$.
 \end{enumerate}
\end{thm}

One may ask if the converse of Theorem \ref{T:303} holds. Before
the end of this section, we give the following Question.
\begin{ques}
 Suppose that $\H=(1,3,h_2,\ldots,h_s)$ is the $h$-vector of a
 {\em level} algebra $A=R/I$ where $R=k[x_1,x_2,x_3]$. Is there a
 level algebra $A$ with the WLP such that $\H$ is the Hilbert function
 of $A$ if $\H=(1,3,h_2,\ldots,h_s)$ satisfies the condition (a),
 (b) and (c) in Theorem \ref{T:303}?

\end{ques}

\section{The Lex-segment Ideals and Graded Non-level Artinian Algebras} \label{Sec:004}

In this section, we shall find an answer to Question~\ref{Q:111}.

\begin{thm}\label{T:403}
Let $R=k[x_1,x_2,x_3]$ and let $\bh=(h_0,h_1,\ldots,h_s)$ be the
$h$-vector of a graded Artinian algebra $A=R/I$ with socle degree $s$. If
$$
h_{d-1}>h_d \qquad \text{and} \qquad h_d=h_{d+1}\leq 2d+3,
$$
then $\H$ is {\bfseries not} level.
\end{thm}

Before we prove this theorem, we consider the following lemmas and the theorems.

\begin{lem} \label{L:403}
Let $J$ be a lex-segment ideal in $R=k[x_1,x_2,x_3]$ such that
$$
\H(R/J,i)=h_i
$$
for every $i\ge 0$. Then
\begin{equation}\label{EQ:403}
\dim_k\left ((J:x_3)/J\right )_{i}=h_{i}-h_{i+1}+(h_{i+1})^{-}
\end{equation}
for such an $i$.
\end{lem}

\begin{proof} First of all, we consider the following exact sequence:
\begin{equation}\label{EQ:311}
0 \rightarrow \left ((J:x_3)/J \right)_{i} \rightarrow
(R/J)_{i} \stackrel{\times x_3}{\longrightarrow} (R/J)_{i+1}
\rightarrow R/(J + (x_3))_{i+1} \rightarrow 0.
\end{equation}
Using equations~\eqref{EQ:310} and \eqref{EQ:311}, we see that
\begin{equation}\label{EQ:312}
\dim_k\left ((J:x_3)/J\right )_{i}=h_{i}-h_{i+1}+(h_{i+1})^{-}
\end{equation}
for every $i\ge 0$ as we desired.
\end{proof}

Since the following lemma is obtained easily from the property of the lex-segment ideal, we shall omit the proof here.

\begin{lem} \label{L:013}
Let $I$ be the  lex-segment ideal in $R=k[x_1,x_2,x_3]$ with Hilbert function $\H=(h_0,h_1,\dots,h_s)$ where  $h_d=d+i$ and $1\le i\le \frac{d^2+d}{2}$.  Then the last monomial of $I_d$ is
$$
\begin{array}{clrlllllllllll}
x_1x_2^{i-1}x_3^{d-i}, & \text{ for } & 1 & \le & i & \le & d, \\
x_1^2x_2^{i-(d+1)}x_3^{(2d-1)-i}, & \text{ for } & d+1 & \le & i & \le & 2d-1,\\
\vdots\\
x_1^{d-1}x_2^{i-\frac{d^2+d-4}{2}}x_3^{\frac{d^2+d-2}{2}-i}, & \text{ for } & \frac{d^2+d-4}{2} & \le & i & \le & \frac{d^2+d-2}{2},\\
x_1^{d}, & \text{ for } &   &   & i & = & \frac{d^2+d}{2}.
\end{array}
$$
\end{lem}

\begin{thm} \label{T:014}
Let $R=k[x_1,x_2,x_3]$ and let $\H=(h_0,h_1,\dots,h_s)$ be the $h$-vector of an Artinian algebra with socle degree $s$ and
$$
h_d=h_{d+1}=d+i, \quad h_{d-1}>h_d, \quad \text{and} \quad j:=h_{d-1}-h_d
$$
for $i=1,2,\dots,\frac{d^2+d}{2}$. Then, for every $1\le k\le d$ and $1\le \ell \le d$,
$$
\begin{array}{llllllllll}
\ds \beta_{1,d+2}
= \begin{cases}
2k-1, & \text{for } \quad
(k-1)d-\frac{k(k-3)}{2} \le i \le (k-1)d-\frac{k(k-3)}{2}+(k-1),
\\
2k, & \text{for } \quad (k-1)d-\frac{k(k-3)}{2}+k\le i \le kd-\frac{(k-1)k}{2}.
\end{cases}\\[5ex]
{\ds \beta_{2,d+2}}
= j+\ell, \hskip 7.9 mm
\text{ for } \quad (\ell-1) d-\frac{(\ell-2)(\ell-1)}{2}<i\le  \ell d-\frac{(\ell-1)\ell}{2}.
\end{array}
$$
\end{thm}

\begin{proof}
Since  $h_d=d+i$, the monomials not in $I_d$ are the last $d+i$ monomials of $R_d$.
By Lemma~\ref{L:013}, the last monomial of $R_1I_d$ is
$$
\begin{array}{clllllllllllll}
x_1x_2^{i-1}x_3^{d-i+1}, & \text{ for } & i=1,\dots, d, \\
x_1^2x_2^{i-(d+1)}x_3^{2d-i}, & \text{ for }  & i=d+1,\dots,2d-1,\\
\vdots\\
x_1^{d-1}x_2^{i-\frac{d^2+d-4}{2}}x_3^{\frac{d^2+d}{2}-i}, & \text{ for } &  i=\frac{d^2+d-4}{2},\ \frac{d^2+d-2}{2}, \\
x_1^{d}x_3, & \text{ for } &    i = \frac{d^2+d}{2}.
\end{array}
$$
In what follows, the first monomial of $I_{d+1}-R_1I_d$ is
\begin{equation}\label{EQ:001}
\begin{array}{cllllllllllll}
x_2^{d+1}, & \text{ for } & i=1, \\
x_1x_2^{i-2}x_3^{(d+2)-i}, & \text{ for } & i=2,\dots,d, \\
\vdots\\
x_1^{d-1}x_2x_3, & \text{ for } & i=\frac{d^2+d-2}{2}, \\
x_1^{d-1}x_2^2, & \text{ for } &   i=\frac{d^2+d}{2}.
\end{array}
\end{equation}
Note that
\begin{equation}\label{EQ:002}
\begin{array}{lllllllllllll}
(d+i)^{\langle d\rangle}
& = & (d+i)+k, & \text{ for }\ i=(k-1)d-\frac{k(k-3)}{2},\dots,
kd-\frac{k(k-1)}{2},\\
&   &          & \text{ and }k=1,\dots,d.
\end{array}
\end{equation}

We now calculate the Betti number
$$
\ds \beta_{1,d+2}
= \ds\sum_{T\in \GG(I)_{d+1}} \binom{m(T)-1}{1}.
$$
Based on equation~\eqref{EQ:001}, we shall find this Betti number of each two cases for $i$ as follows.
\begin{enumerate}

\item[]{\em Case 1-1.}  $i=(k-1)d-\frac{k(k-3)}{2}$ and $k=1,2,\dots,d$.

By equation~\eqref{EQ:002}, $I_{d+1}$ has $k$-generators, which are
$$
x_1^{k-1} x_2^{(d+2)-k},x_1^{k-1} x_2^{(d+1)-k}x_3,\dots,x_1^{k-1} x_2^{(d+3)-2k} x_3^{k-1}.
$$
By the similar argument, $I_{d+1}$ has $k$-generators including the element $x_1^{k-1} x_2^{(d+2)-k}$ for $i=(k-1)d-\frac{k(k-3)}{2}+1,\dots, (k-1)d-\frac{k(k-3)}{2}+(k-1)$. Hence we have that
$$
\ds \beta_{1,d+2}
= \ds\sum_{T\in \GG(I)_{d+1}} \binom{m(T)-1}{1}=2\times (k-1)+1=2k-1.
$$

\item[]{\em Case 1-2.} $i=(k-1)d-\frac{k(k-3)}{2}+k=(k-1)d-\frac{k(k-5)}{2},\dots,
kd-\frac{k(k-1)}{2}$ and
$k=1,2,\dots,d$.

By equation~\eqref{EQ:002}, $I_{d+1}$ has $k$-generators, which are
$$
x_1^{k} x_2^{i-\left((k-1)d-\frac{k^2-3k-2}{2}\right)}
x_3^{kd-\frac{k^2-k-4}{2}-i},\dots,
x_1^{k} x_2^{i-\left((k-1)d-\frac{k(k-5)}{2}\right)}
x_3^{\left(kd-\frac{k(k-3)}{2}+1\right)-i}.
$$
Hence we have that
$$
\ds \beta_{1,d+2}
= \ds\sum_{T\in \GG(I)_{d+1}} \binom{m(T)-1}{1}=2\times k=2k.
$$
\end{enumerate}

Now we move on to the Betti number:
$$
\ds \beta_{2,d+2}=\sum_{T\in \GG(I)_d} \binom{m(T)-1}{2}.
$$
Recall $h_d=d+i$ and $j:=h_{d-1}-h_d$. The computation of the Betti number of this case is much more complicated, and thus we shall find the Betti number of each four cases based on $i$ and $j$.

\begin{enumerate}

\item[]{\em Case 2-1.} $(\ell-1) d-\frac{(\ell-2)(\ell-1)}{2}<i< \ell d-\frac{(\ell-1)\ell}{2}$ and $\ell=1,2,\dots,d$.

The last monomial of $I_d$ for this case is
$$
x_1^{\ell} x_2^{i-(\ell-1) d+\frac{\ell(\ell-3)}{2}} x_3^{\ell d -\frac{(\ell-1)\ell}{2}-i}.
$$

\begin{enumerate}

\item[]{\em Case 2-1-1.}  $(k-1)d-\frac{(k-1)k}{2}<i+j<kd-\frac{k(k+1)}{2}$ and $k=\ell,\ell+1,\dots,d$.

Since first monomial of $I_d-R_1I_{d-1}$ is
$$
x_1^{k} x_2^{(i+j)-\left((k-1)d-\frac{(k-2)(k+1)}{2}\right)}x_3^{\left(kd-\frac{(k-1)(k+2)}{2}\right)-(i+j)},
$$
we have $(j+k)$-generators in $I_d$ as follows:
$$
\begin{array}{llllllllllllllll}
x_1^{k} x_2^{(i+j)-\left((k-1)d-\frac{(k-2)(k+1)}{2}\right)}x_3^{\left(kd-\frac{(k-1)(k+2)}{2}\right)-(i+j)}, \dots,x_1^{k}x_3^{d-k},\\
x_1^{(k-1)} x_2^{d-(k-1)}, x_1^{(k-1)} x_2^{(d-1)-(k-1)}x_3,\dots, x_1^{(k-1)} x_3^{d-(k-1)}, \\
\hskip 2.5 true cm \vdots \\

x_1^{\ell+1} x_2^{(d-1)-\ell}, x_1^{\ell+1} x_2^{(d-2)-\ell}x_3,
\dots, x_1^{\ell+1} x_3^{(d-1)-\ell}\\
x_1^{\ell} x_2^{d-\ell},  \dots,
x_1^{\ell} x_2^{i-(\ell-1) d+\frac{\ell(\ell-3)}{2}} x_3^{\ell d -\frac{(\ell-1)\ell}{2}-i}\\
\end{array}
$$
and thus
$$
\ds \beta_{2,d+2}=\sum_{T\in \GG(I)_d} \binom{m(T)-1}{2}=j+\ell.
$$

\item[]{\em Case 2-1-2.} $i+j=(k-1)d-\frac{(k-1)k}{2}$ and $k=\ell+1,\dots,d$.

The first monomial of $I_d-R_1I_{d-1}$ is
$$
x_1^{k-1} x_2^{d-(k-1)},
$$
and hence we have $(j+k)$-generators in $I_d$ as follows:
$$
\begin{array}{llllllllllllllll}
x_1^{k-1} x_2^{d-(k-1)}, x_1^{k-1} x_2^{(d-1)-(k-1)}x_3,\dots, x_1^{k-1} x_3^{d-(k-1)}, \\
\hskip 2.5 true cm \vdots \\
x_1^{\ell+1} x_2^{(d-1)-\ell}, x_1^{\ell+1} x_2^{(d-2)-\ell}x_3,
\dots, x_1^{\ell+1} x_3^{(d-1)-\ell}\\
x_1^{\ell} x_2^{d-\ell},  \dots,
x_1^{\ell} x_2^{i-(\ell-1) d+\frac{\ell(\ell-3)}{2}} x_3^{\ell d -\frac{(\ell-1)\ell}{2}-i}\\
\end{array}
$$
and thus
$$
\ds \beta_{2,d+2}=\sum_{T\in \GG(I)_d} \binom{m(T)-1}{2}=j+\ell.
$$

\end{enumerate}

\item[]{\em Case 2-2.} $i=\ell d-\frac{(\ell-1)\ell}{2}$ and $\ell=1,2,\dots,d$.

The last monomial of $I_d$ is
$$
x_1^\ell x_2^{d-\ell}.
$$

\begin{enumerate}

\item[]{\em Case 2-2-1.} $(k-1)d-\frac{(k-1)k}{2}<i+j<kd-\frac{k(k+1)}{2}$ and $k=\ell+1,\dots,d$.

Since the first monomial of $I_d-R_1I_{d-1}$ is
$$
x_1^{k} x_2^{(i+j)-\left((k-1)d-\frac{(k-2)(k+1)}{2}\right)}x_3^{\left(kd-\frac{(k-1)(k+2)}{2}\right)-(i+j)},
$$
we have $(j+k)$-generators in $I_d$ as follows:
$$
\begin{array}{llllllllllllllll}
x_1^{k} x_2^{(i+j)-\left((k-1)d-\frac{(k-2)(k+1)}{2}\right)}x_3^{\left(kd-\frac{(k-1)(k+2)}{2}\right)-(i+j)}, \dots,x_1^{k}x_3^{d-k},\\
x_1^{(k-1)} x_2^{d-(k-1)}, x_1^{(k-1)} x_2^{(d-1)-(k-1)}x_3,\dots, x_1^{(k-1)} x_3^{d-(k-1)}, \\
\hskip 2.5 true cm \vdots \\
x_1^{\ell+1} x_2^{(d-1)-\ell}, x_1^{\ell+1} x_2^{(d-2)-\ell}x_3,
\dots, x_1^{\ell+1} x_3^{(d-1)-\ell}\\
x_1^\ell x_2^{d-\ell},
\end{array}
$$
and thus
$$
\ds \beta_{2,d+2}=\sum_{T\in \GG(I)_d} \binom{m(T)-1}{2}=j+\ell.
$$

\item[]{\em Case 2-2-2.} $i+j=(k-1)d-\frac{(k-1)k}{2}$ and $k=\ell+1,\dots,d$.

The first monomial of $I_d-R_1I_{d-1}$ is
$$
x_1^{(k-1)} x_2^{d-(k-1)},
$$
and hence we have $(j+k)$-generators in $I_d$ as follows:
$$
\begin{array}{llllllllllllllll}
x_1^{(k-1)} x_2^{d-(k-1)}, x_1^{(k-1)} x_2^{(d-1)-(k-1)}x_3,\dots, x_1^{(k-1)} x_3^{d-(k-1)}, \\
\hskip 2.5 true cm \vdots \\
x_1^{\ell+1} x_2^{(d-1)-\ell}, x_1^{\ell+1} x_2^{(d-2)-\ell}x_3,
\dots, x_1^{\ell+1} x_3^{(d-1)-\ell}\\
x_1^\ell x_2^{d-\ell},
\end{array}
$$

and thus
$$
\ds \beta_{2,d+2}=\sum_{T\in \GG(I)_d} \binom{m(T)-1}{2}=j+\ell,
$$
\end{enumerate}
\end{enumerate}
as we wished.
\end{proof}

\begin{thm} \label{T:407}
Let $\H$ be as in equation~\eqref{EQ:011} and $A=R/I$ be an algebra with Hilbert function $\H$ such that $\beta_{1,d+2}(I^{\rm lex})=\beta_{2,d+2}(I^{\rm lex})$ for some $d<s$. Then $A$ is {\bfseries not} level.
\end{thm}

\begin{proof} Let $L$ be a general linear form of $A$. By Lemma \ref{lemma 1} (b), note that if $d\ge r_1(A)$, then
$$
\dim_k(0:L)_{d-1}\geq h_{d-1}-h_d>0 \quad \text{ and } \quad \dim_k(0:L)_{d}=h_d-h_{d+1}=0,
$$
and thus,   by Lemma \ref{lemma 1} (a), $R/I$ is not level. Hence we
assume that $d<r_1(A)$ and $A$ is a graded level algebra
having Hilbert function $\H$.
 Let $\bar{I}=(I_{\leq d+1})$.

\medskip

\noindent{\em Claim.} $\beta_{1,d+3}(\Gin(\bar{I}))=0$ and
$\beta_{2,d+3}(\Gin(\bar{I}))>0$.

\medskip

\ni
{\em Proof of Claim.}
First we shall show that $\beta_{1,d+3}(\Gin(\bar{I}))=0$. By Lemma~\ref{L:209} (a),
$$
\beta_{1,d+2}({I}^{\rm lex} )=\beta_{2,d+2}({I}^{\rm lex} ),
$$
and we have that
\begin{equation} \label{EQ:413}
\begin{array}{rrlllllllllllll}
&
\beta_{1,d+2}(I^{\rm lex})-\beta_{1,d+2}(I)
& = & [\beta_{0,d+2}(I^{\rm lex})-\beta_{0,d+2}(I)]
+[\beta_{2,d+2}(I^{\rm lex})-\beta_{2,d+2}(I)]\\
\Rightarrow &
 -\beta_{1,d+2}(I)
& = & [\beta_{0,d+2}(I^{\rm lex})-\beta_{0,d+2}(I)]
 -\beta_{2,d+2}(I).
\end{array}
\end{equation}
Moreover, since $A=R/I$ is level, we know that $\beta_{2,d+2}(I)=0$, and hence rewrite equation~\eqref{EQ:413} as
$$
0 \le [\beta_{0,d+2}(I^{\rm lex})-\beta_{0,d+2}(I)]= -\beta_{1,d+2}(I) \le 0,
$$
which follows from Lemma~\ref{L:208} (b) that
$$
\beta_{0,d+2}(I^{\rm lex})-\beta_{0,d+2}(I)=
\beta_{0,d+2}(\bar{I}^{\rm lex})=0.
$$
Also, by Lemma~\ref{L:208} (a), we have
$$
\beta_{0,d+2}({\rm Gin}(\bar I))\le \beta_{0,d+2}(\bar{I}^{\rm lex}))=0,
\quad \text{i.e.,} \quad \beta_{0,d+2}({\rm Gin}(\bar I))=0.
$$
Since ${\rm Gin}(\bar I)$ is a Borel fixed monomial ideal, by Theorem~\ref{EK},
$$
\beta_{1,d+3}({\rm Gin}(\bar I))=0.
$$

\medskip

Now we shall prove that $\beta_{2,d+3}(\Gin(\bar{I}))>0$. Let $J={\rm Gin}(\bar I)$. Consider the following exact sequence
$$
0 \rightarrow \left ((J:x_3)/J \right )_{d} \rightarrow (R/J)_{d}
\stackrel{\times x_3}{\longrightarrow} (R/J)_{d+1} \rightarrow
(R/J + (x_3))_{d+1} \rightarrow 0.
$$
Since $d< r_1(A)$, we know that
$$
\begin{array}{llllllllllllllllllllllllll}
\dim_k\left ((J:x_n)/J \right )_{d}
& = &\,h_d-h_{d+1}+ \dim_k((R/J + (x_3))_{d+1})\\[1ex]
& = & \dim_k((R/J + (x_n))_{d+1}) \qquad (\because h_d=h_{d+1}) \\[1ex]
& \ne&  0.
\end{array}
$$
By Lemma~\ref{strongly stable},
$$
\mathcal G(J)_{d+1}=\mathcal G({\rm Gin}(\bar I))_{d+1}\ne \varnothing,
$$
and so there is a monomial
$T\in \mathcal G({\rm Gin}(\bar I))_{d+1}$ such that $x_3\,|\, T$. In other words,
$$
\beta_{2,d+3}(\Gin(\bar{I}))>0,
$$
as we desired.

\medskip

By the above claim and a cancellation principle, $R/{\bar I}$ has a socle element in degree $d$,  and thus $R/I$ has such a socle element in degree $d$ since $R/I$ and $R/{\bar I}$ agree in degrees $\le d+1$, and hence $A$ cannot be level, as we wished.
\end{proof}

Now we are ready to prove Theorem~\ref{T:403}.

\begin{proof}[Proof of Theorem~\ref{T:403}]
Let $\H$ and $j$ be as in Theorem~\ref{T:014} and let $h_{d}=d+i$ for $-(d-1)\le i \le d+3$.

By Proposition 3.8 in \cite{GHMS}, this theorem holds for $-(d-1) \le i \le 1$. It suffices, therefore, to prove this theorem for $2 \le i \le d+3$. By Theorem~\ref{T:014}, we have
\begin{equation} \label{EQ:003}
\begin{array}{lllllllllllll}
\beta_{1,d+2}(I^{\rm lex})=
\begin{cases}
2, & \text{for }  i=2,\dots,d, \\
3, & \text{for }  i=d+1,d+2, \\
4, & \text{for }  i=d+3,
\end{cases}
\quad \text{and} \\[5ex]
\beta_{2,d+2}(I^{\rm lex})=
\begin{cases}
j+1, & \text{for } i=2,\dots,d, \\
j+2, & \text{for }  i=d+1, d+2, d+3.
\end{cases}
\end{array}
\end{equation}

Note that if either $j\ge 3$ and $2\le i\le d+3$ or $j=2$ and $2\le i \le d+2$, then $\H$ is not level since $\beta_{2,d+2}(I^{\rm lex})>\beta_{1,d+2}(I^{\rm lex})$.

\medskip

Now suppose either $j=1$ and $2\le i\le d+2$ or $j=2$ and $i=d+3$. By equation~\eqref{EQ:003}, we have
$$
\beta_{1,d+2}(I^{\rm lex})=\beta_{2,d+2}(I^{\rm lex})=
\begin{cases}
2, & \text{for } j=1 \text{ and } i=2,\dots,d, \\
3, & \text{for } j=1 \text{ and } i=d+1,d+2, \\
4, & \text{for } j=2 \text{ and } i=d+3.
\end{cases}
$$
Thus, by Theorem~\ref{T:407}, $\H$ cannot be level.

\medskip

It is enough, therefore,  to show the case $j=1$ and $i=d+3$.  Assume there exists a level algebra $R/I$ with Hilbert function $\H$.
Applying equation~\eqref{EQ:003} again, we have
\begin{equation} \label{EQ:409}
\beta_{1,d+2}(I^{\rm lex})=\beta_{2,d+2}(I^{\rm lex})+1=4.
\end{equation}

Note $h_{d-1}=2d+4$ and $h_d=h_{d+1}=2d+3$ in this case.  By equation~\eqref{EQ:409}, the Betti diagram of $R/I^{\rm lex}$ is as follows
\begin{table}[ht]
\begin{center}
\begin{verbatim}
                              total: 1   -   -   -
                             ------------------------
                                  0: 1   -   -   -
                                  1: -   -   -   -
                                          ......
                                d-1: -   *   *   3
                                  d: -   *   4   *
                                d+1: -   *   *   *
                                          ......
\end{verbatim}
\vskip .5pc
\caption{Betti diagram of $R/I^{\rm lex}$} \label{TB:002}
\vskip -1pc
\end{center}
\end{table}

Moreover, by Lemmas~\ref{strongly stable} and~\ref{L:403},
\begin{equation} \label{EQ:410}
\begin{array}{llllllllllllllllll}
\dim_k ((I^{\rm lex}:x_3)/I^{\rm lex})_{d}
& = & \Big|\left\{\,T\in \GG(I^{\rm lex})_{d+1} \,\Big|\, x_3\, \big|\, T\,\right\}\Big| \\[2ex]
& = & h_{d}-h_{d+1}+(h_{d+1})^{-} \\[1ex]
& = & (h_{d+1})^{-} \\[1ex]
& = & \ds\left(\binom{d+2}{d+1}+\binom{d+1}{d}\right)^{-} \\[2ex]
& = & 2.
\end{array}
\end{equation}

Hence, using equation~\eqref{EQ:410},  we can rewrite Table~\ref{TB:002} as
\begin{table}[ht]
\begin{center}
\begin{verbatim}
                              total: 1   -   -   -
                             ------------------------
                                  0: 1   -   -   -
                                  1: -   -   -   -
                                          ......
                                d-1: -   *   *   3
                                  d: -   2   4   2
                                d+1: -   *   *   *
                                          ......
\end{verbatim}
\vskip .5pc
\caption{Betti diagram of $R/I^{\rm lex}$} \label{TB:003}
\vskip -1pc
\end{center}
\end{table}

Let $J:=(I_{\le d+1})^{\rm lex}$. Note  $I^{\rm lex}$ and $J$ agree in degree $\le d+1$.  Hence we can write the Betti diagram of $R/J$ as


\begin{table}[ht]
\begin{center}
\begin{verbatim}
                              total: 1   -   -   -
                             ------------------------
                                  0: 1   -   -   -
                                  1: -   -   -   -
                                          ......
                                d-1: -   *   *   3
                                  d: -   2   4   2
                                d+1: -   a   b   *
                                          ......
\end{verbatim}
\vskip .5pc
\caption{Betti diagram of $R/J$} \label{TB:004}
\vskip -1pc
\end{center}
\end{table}

Since $R/I$ is level and $(I_{\le d+1})$ has no generators in degree $d+2$, we have  $a=0$ or $1$.

\noindent
{\em Case 1.} Let $a=0$. Then, by Theorem~\ref{EK}, we have $b=0$. Since $J$ and $(I_{\le d+1})$ agree in degree $\le d+1$,
$$
\beta_{2,d+3}(J)=\beta_{2,d+3}((I_{\le d+1}))=2.
$$
This means  $R/(I_{\le d+1})$ has two dimensional socle elements in degree $d$, so does $R/I$, which is a contradiction.

\medskip

\noindent
{\em Case 2.} Let $a=1$, then $J$ has one generator in degree $d+2$.
By Lemmas~\ref{strongly stable} and~\ref{L:403},
\begin{equation}\label{EQ:411}
\begin{array}{llllllllllllllllllll}
\dim_k ((J:x_3)/J)_{d+1}
& = & \Big|\left\{\, T \in \GG(J)_{d+2} \,\big|\, x_3 \,\big|\, T\, \right\}\Big|
\\[1ex]
& = & h_{d+1}-h_{d+2}+(h_{d+2})^{-}
\end{array}
\end{equation}
where $h_{d+2}=\H(R/J,d+2)=h_{d+1}^{\langle d+1\rangle}-1=(2d+3)^{\langle d+1\rangle}-1=2d+4$. Hence we obtain $(h_{d+2})^{-}=(2d+4)^{-}=1$, and by equation~\eqref{EQ:411}
$$
\dim_k ((J:x_3)/J)_{d+1}=0.
$$

\medskip

Applying Theorem~\ref{EK} again, we find
$$
b=\beta_{1,d+3}(J)=\sum_{T\in \GG(J)_{d+2}}\binom{m(T)-1}{1}=1
$$
since $x_1^{d+2}\notin   \GG(J)_{d+2}$. Thus $R/J$  has at least one socle element in degree $d$, and so does $R/(I_{\le d+1})$. Since $R/I$ and $R/(I_{\le d+1})$ agree in degree $\le d+1$, $R/I$ has such a socle element, a contradiction, which completes the proof.
\end{proof}

The following example shows a case where $j=1$ and $h_d=2d+3$ in Theorem~\ref{T:403}.

\begin{exmp} \label{R:058} Let $I$ be the lex-segment ideal in $R=k[x_1,x_2,x_3]$ with Hilbert function
$$
\begin{array}{lllllllllllllllllllll}
\H & : & 1 & 3 & 6 & 10 & 15 & 21 & 18 & 17 & 17 & 0 & \rightarrow \ .
\end{array}
$$
Note that $h_7=17=2\times 7+3=2d+3$, which satisfies the condition in Theorem~\ref{T:403}, and $j=h_6-h_7=18-17=1$. Hence
any Artinian algebra having Hilbert function $\H$ cannot be level.
\end{exmp}

Inverse systems can also be used to produce new level algebras from known level algebras. This method is based on the idea of {\em Macaulay's Inverse Systems} (see \cite{Ger} and \cite{I-K} for details). We want to recall some results from \cite{Ia:1}. Actually, Iarrobino shows an even stronger result and the application to level algebras is:

\begin{thm}[Theorem 4.8A, \cite{Ia:1}]  \label{T:405}
Let $R=k[x_1,\dots,x_r]$ and $\H'=(h_0,h_1,\dots,h_e)$ be the $h$-vector of a level algebra $A=R/{\rm Ann}(M)$. Then, if $F$ is a generic form of degree $e$, the level algebra $R/{\rm Ann}(\langle M,F\rangle)$ has $h$-vector $\H=(H_0,H_1,\dots,H_e)$, where, for $i=1,\dots,e$,
$$
H_i=\min \left \{h_i+\binom{(r-1)+(e-i)}{(e-i)}, \binom{(r-1)+i}{i}\right\}.
$$
\end{thm}

The following example is another case of a level O-sequence of codimension $3$ of type in equation~\eqref{EQ:011} satisfying $h_d=2d+4$.

\begin{exmp} \label{EX:410}
Consider a level O-sequence $(1, 3, 5, 7, 9, 11, 13)$ of codimension $3$. By Theorem~\ref{T:405}, we obtain the following level O-sequence:
$$
(1, 3, 6, 10, 15, 14, 14).
$$
Then $14=2\times 5+4$, which shows  there exists a level O-sequence of codimension $3$ of type in equation~\eqref{EQ:011} when $h_d=2d+4$.
\end{exmp}

In general, we can construct a level O-sequence of codimension $3$ of type in equation~\eqref{EQ:011} satisfying $h_d=2d+4$ for every $d\ge 5$ as follows.

\begin{pro} \label{P:411}
There exists a level O-sequence of codimension $3$ of type in equation~\eqref{EQ:011} satisfying $h_d=2d+4$ for every $d\ge 5$.
\end{pro}

\begin{proof} Note that, from Example~\ref{EX:410}, this proposition holds for $d=5$.

Now assume  $d\ge 6$. Consider a level O-sequence $h=(1, 3, 5, 7,\dots,\overset{d\text{-th}}{2d+1},\overset{(d+1)\text{-st}}{2d+3})$ where $d\ge 6$. Since
$$
\begin{array}{llllllllllllllllllllllll}
&   &  \left(h_{i}+\binom{d+3-i}{d+1-i}\right)-\binom{i+2}{i}\\[2ex]
& = &  \left(2i + 1+\frac{(d + 3 - i)(d + 2 - i)}{2}\right)-\frac{(i+1)(i+2)}{2} \\[2ex]
& = &  \frac{(2+d)(3+d-2i)}{2}\ge 0,
\end{array}
$$
for every $i=0,1,\dots,d-3$, we have
$$
\begin{array}{llllllllllllllllllllllllllllll}
H_{i} & = & \min\left\{h_{i}+\binom{d+3-i}{d+1-i},\binom{i+2}{i}\right\}\\[2ex]
        & = & \min\left\{2i + 1+\frac{(d + 3 - i)(d + 2 - i)}{2}, \frac{(i+1)(i+2)}{2}\right\} \\[2ex]
& = & \frac{(i+1)(i+2)}{2}.
\end{array}
$$
Hence, by Theorem~\ref{T:405},
 we obtain a level O-sequence $\H=(H_0,H_1,\dots,H_d,H_{d+1})$ as follows:
$$
\begin{array}{lllllllllllllllllllllllllllll}
H_0 & = & 1, \\
H_1 & = & 3, \\
    & \vdots & \\
H_i & = & \frac{(i+1)(i+2)}{2}, \\
    & \vdots & \\
H_{d-2} & = & \min\left\{h_{d-2}+\binom{5}{3},\binom{d}{d-2}\right\}
        & = & \min\left\{2d+7, \frac{(d-1)d}{2}\right\} & = & 2d+7,\\[1ex]
H_{d-1} & = & \min\left\{h_{d-1}+\binom{4}{2},\binom{d+1}{d-1}\right\}
        & = & \min\left\{2d+5, \frac{d(d+1)}{2}\right\} & = & 2d+5,\\[1ex]
H_{d} & = & \min\left\{h_{d}+\binom{3}{1},\binom{d+2}{d}\right\}
        & = & \min\left\{2d+4, \frac{(d+1)(d+2)}{2}\right\} & = & 2d+4,\\[1ex]
H_{d+1} & = & \min\left\{h_{d+1}+\binom{2}{0},\binom{d+3}{d+1}\right\}
        & = & \min\left\{2d+4, \frac{(d+2)(d+3)}{2}\right\} & = & 2d+4,
\end{array}
$$
as we wished.
\end{proof}

\begin{rem} \label{R:412}
By the same idea as in the proof of Proposition~\ref{P:411}, we can construct a level O-sequence of codimension $3$ of type in equation~\eqref{EQ:011} satisfying
$$
2d+(k+1)=H_{d-1}>H_d=H_{d+1}= 2d+k, \quad (\begin{matrix} 5\le k \le \frac{d^2-3d+2}{2}\end{matrix}).
$$
For example, if we use
$$
h=(1, 3, 6, \dots,\overset{(d-1)\text{-st}}{2d+(k-5)},\overset{d\text{-th}}{2d+(k-3)},\overset{(d+1)\text{-st}}{2d+(k-1)}),
$$
then we construct a level O-sequence of codimension $3$ of type in equation~\eqref{EQ:011} satisfying
$$
\begin{array}{lllllllllllllllllllll}
H_{d-1}
& = & \min\left\{h_{d-1}+\binom{4}{2},\binom{d+1}{d-1}\right\}
& = & \min\left\{2d+(k+1), \frac{d(d+1)}{2}\right\} & = & 2d+(k+1), \\
&   & (\because k \le \frac{d^2-3d-2}{2}), \\[1ex]
H_{d}
& = & \min\left\{h_{d}+\binom{3}{1},\binom{d+2}{d}\right\}
& = & \min\left\{2d+k, \frac{(d+1)(d+2)}{2}\right\} & = & 2d+k, \\[1ex]
H_{d+1}
& = & \min\left\{h_{d+1}+\binom{2}{0},\binom{d+3}{d+1}\right\}
& = & \min\left\{2d+k, \frac{(d+2)(d+3)}{2}\right\} & = & 2d+k,
\end{array}
$$
as we desired.
\end{rem}

Using Theorem~\ref{T:403}, we know that some non-unimodal O-sequence of codimension $3$ cannot be level as follows.

\begin{cor} \label{C:510}
Let $\H=\{h_i\}_{i\ge 0}$ be an O-sequence with $h_1=3$. If
$$
h_{d-1}>h_d, \quad h_d\le 2d+3, \quad \text{ and } \quad h_{d+1}\ge h_d
$$
for some degree $d$, then $\H$ is not level.
\end{cor}

\begin{proof} Note that, by the proof of Theorem~\ref{T:403}, any graded ring with Hilbert function
$$
\begin{matrix}
\H' & : & h_0 & h_1 & \cdots & h_{d-1} & h_d & h_d & \rightarrow
\end{matrix}
$$
has a socle element in degree $d-1$.

Now let $A=\bigoplus_{i\ge 0} A_i$ be a graded ring with Hilbert function $\H$.  If $A_{d+1}=\langle f_1,f_2,\dots,f_{h_{d+1}}\rangle$ and $I=(f_{h_d+1},\dots,f_{h_{d+1}})\bigoplus_{j\ge d+2} A_j$, then a graded ring $B=A/I$ has Hilbert function
$$
\begin{matrix}
h_0 & h_1 & \cdots & h_{d-1} & h_d & h_d.
\end{matrix}
$$
Hence $B$ has a socle element in degree $d-1$ or $d$ by Theorem~\ref{T:403}. Since $A_i=B_i$ for every $i\le d$, $A$ also has the same socle element in degree $d-1$ or $d$ as $B$, and thus $\H$ is not level as we wished.
\end{proof}

The following is an example of a non-level and non-unimodal O-sequence of codimension $3$ satisfying the condition of Corollary~\ref{C:510}.

\begin{exmp} Consider an O-sequence
$$
\begin{array}{lllllllllllllllllllll}
\H & : & 1 & 3 & 6 & 10 & 15 & 20 & 18 & 17 & h_8 & \cdots .
\end{array}
$$
There are only $3$ possible O-sequences such that $h_8\ge h_7=17$ since $h_8\le h_7^{\langle 7\rangle}=17^{\langle 7\rangle}=19$. By Theorem~\ref{T:403}, $\H$ is not level if $h_8=h_7=17$. The other two non-unimodal O-sequences, by Corollary~\ref{C:510},
$$
\begin{array}{lllllllllllllllllllll}
1 & 3 & 6 & 10 & 15 & 20 & 18 & 17 & 18 & \cdots \quad \text{and}\\
1 & 3 & 6 & 10 & 15 & 20 & 18 & 17 & 19 & \cdots
\end{array}
$$
cannot be level either.
\end{exmp}

\bibliographystyle{amsalpha}

\end{document}